\documentclass[a4paper,twoside]{article}
\usepackage{amssymb,amsmath,amscd,amsthm}

\input{epsf}

\newcommand{\str}{^\prime}
\newcommand{\res}[1]{\!\mid_{#1}}
\newcommand{\card}[1]{{\mid\! #1 \!\mid}}
\newcommand{\pro}[2]{\langle #1, #2 \rangle}
\newcommand{\qu}[1]{{\overline{#1}}}

\def\lolra{\Longleftrightarrow}

\def\N{{\mathbb N}}
\def\P{{\mathbb P}}
\def\Q{{\mathbb Q}}
\def\R{{\mathbb R}}
\def\Z{{\mathbb Z}}

\def\F{{\cal F}}

\def\calN{{\cal N}}

\def\TT{{\mathbb T}}
\def\T{{\cal T}}

\def\V{{{\cal V}}}
\def\W{{\cal W}}

\def\conv{{\rm conv}}
\def\aff{{\rm aff}}
\def\lin{{\rm lin}}
\def\pos{{\rm pos}}

\def\relint{{\rm relint}}
\def\rank{{\rm rank}}

\def\Hom{{\textup{Hom}}}

\def\DivTX{{\rm Div}^{\TT}(X)}

\def\Cl{\textup{Cl}}

\def\NR{{N_{\R}}}
\def\NQ{{N_{\Q}}}
\def\MR{{M_{\R}}}
\def\MQ{{M_{\Q}}}

\def\fan{\triangle}
\def\fans{\fan\str}
\def\taus{{\tau\str}}
\def\fanqu{{\qu{\fan}}}
\def\fanv{{\fan_v}}
\def\fanT{{\fan_\T}}
\def\rand{\partial}
\def\randp{{\rand P}}

\def\randx{{\rand x}}
\def\rays{{\fan(1)}}
\def\rayss{{\fans(1)}}

\def\intr{{\rm int}}
\def\st{{\rm st}}
\def\VP{{\V(P)}}
\def\randPM{{\rand P \cap M}}

\def\XT{{X_\T}}

\def\Xv{{X_v}}

\def\Vt{{\V_\tau}}

\def\vt{{v_\tau}}
\def\vts{{v_\taus}}
\def\tg{{\tau\str}}
\def\vtg{{v_\tg}}

\def\KX{{K_X}}

\def\jX{{j_X}}
\def\disc{{\rm discr}}
\def\us{{u_\sigma}}
\def\Gf{{G_\fan}}
\def\Gs{{G_\sigma}}
\def\Qf{{Q_\fan}}
\def\Pf{{P_\fan}}

\def\Zt{{\cal Z}}
\def\id{\textup{id}}

\newtheorem{theorem}{Theorem}[section]
\newtheorem{definition}[theorem]{Definition}
\newtheorem{remark}[theorem]{Remark}
\newtheorem{example}[theorem]{Example}
\newtheorem{corollary}[theorem]{Corollary}
\newtheorem{proposition}[theorem]{Proposition}
\newtheorem{lemma}[theorem]{Lemma}
\newtheorem{conjecture}[theorem]{Conjecture}
\newtheorem*{acknowledgement}{Acknowledgement}

\title{Gorenstein Toric Fano Varieties}

\author{{\sc Benjamin Nill} \\
\small  {\em Mathematisches Institut, Universit\"at T\"ubingen}   \\
\small  {\em Auf der Morgenstelle 10,  72076  T\"ubingen, Germany}  \\
\small  {\em e-mail: nill@everest.mathematik.uni-tuebingen.de} \\
 }

\begin{document}

\date{}

\maketitle

\begin{abstract}
We investigate Gorenstein toric Fano varieties by combinatorial me\-thods using 
the notion of a reflexive polytope which appeared in connection to mirror symmetry. 
The paper contains generalisations of tools and previously known results for 
nonsingular toric Fano varieties. As applications we obtain new classification 
results, bounds of invariants and formulate conjectures concerning combinatorial and 
geometrical properties of reflexive polytopes.
\end{abstract}

\section*{Introduction}
\label{intro}
Gorenstein toric Fano Varieties, i.e., normal toric varieties whose anticanonical divisor is 
an ample Cartier divisor, correspond to reflexive polytopes introduced by Batyrev in 
\cite{Bat94}. Reflexive polytopes are lattice polytopes 
containing the origin in their interior such that the dual polytope is also a lattice polytope. 
In \cite{Bat94} it was shown that the associated varieties are ambient spaces for Calabi-Yau hypersurfaces and 
together with their duals naturally yield candidates for mirror symmetry pairs. This has raised interest 
in this special class of lattice polytopes among physicists. 

It is known that in fixed dimension $d$ there is only a finite number of isomorphism classes of 
$d$-dimensional reflexive polytopes. 
Using a computer program Kreuzer and Skarke classified all $d$-dimensional reflexive polytopes 
for $d \leq 4$ \cite{KS98,KS00a,KS02}. They found $16$ isomorphism classes of $d$-dimensional 
reflexive polytopes for $d=2$, $4319$ for $d=3$, 
and $473800776$ for $d=4$. This vast work gives mathematicians many examples of reflexive 
polytopes to test conjectures on.

There are many papers devoted to the investigation and classification of nonsingular toric 
Fano varieties \cite{WW82,Bat82,Bat99,Sat00,Deb01,Cas03a,Cas03b}. In this article we present 
new classification results, bounds of invariants and conjectures concerning Gorenstein toric Fano varieties 
by investigating combinatorial and geometrical properties of reflexive polytopes. Thereby we get not only 
results in higher dimensions but also explanations for interesting phenomena observed in lower dimensions.

In the forthcoming paper \cite{Nil04} the results achieved 
will be applied to investigate criteria for the group of automorphisms 
of a Gorenstein toric Fano variety to be reductive. These include the answer to a question 
related to the existence of an Einstein-K\"ahler metric. Furthermore 
a sharp bound on the dimension of the connected component of the 
reductive group of automorphisms will be proven.

The structure of this paper is as follows: 

In section 1 we establish our notations and repeat the basic definitions. 

In sections 2 and 3 two elementary technical tools are investigated and generalised that 
were previously already successfully used to investigate and classify nonsingular toric Fano varieties 
\cite{Bat99,Sat00,Cas03b}.

In section 2 we investigate the properties of the projection of a reflexive polytope along 
a lattice point on the boundary. Thereby we can relate properties of a 
Gorenstein toric Fano variety to that of a lower-dimensional toric Fano variety (e.g., Cor. \ref{corXv}). 
As an application we give a new, purely combinatorial proof (Cor. \ref{nonsi}) of a result due to Batyrev 
\cite[Prop. 2.4.4]{Bat99} saying that the anticanonical 
class of a torus-invariant prime divisor of a nonsingular toric Fano variety is always numerically effective. 

In section 3 we consider pairs of lattice points on the boundary of a reflexive polytope and 
show that in this case there exists a generalisation of the notion of a 
primitive relation as introduced by Batyrev in \cite{Bat91}. 
As an application we show that there are 
constraints on the combinatorics of a reflexive polytope, in particular on 
the diameter of the edge-graph of a  simplicial reflexive polytope (Cor. \ref{comb}). 
Thereby we can prove that certain combinatorial types of polytopes 
cannot be realized as reflexive polytopes (Cor. \ref{six}, Cor. \ref{noiso}). 

In sections 4 and 5 we present applications of the results of sections 2 and 3.

In section 4 we give a short review of classification results of reflexive polytopes in low dimensions. 
We give a concise proof of the classification of reflexive polygons (Prop. \ref{two}), and a scetch of the proof 
of the classification of three-dimensional Gorenstein toric Fano varieties with terminal singularities 
(Thm. \ref{100}).

In section 5 we formulate two conjectures (Conj. \ref{vertconj}, Conj. \ref{mainconj}) 
on the maximal number of vertices of a reflexive polytope, respectively a 
simplicial reflexive polytope. 
The currently best upper bound on the number of vertices of a smooth Fano polytope is due to Debarre 
\cite[Thm. 8]{Deb01}. Here we give a generalisation to the case of a reflexive polytope (Thm. \ref{bounds}). 
The main result of the section is the verification of the conjecture on the maximal number of 
vertices of a simplicial 
reflexive polytope under the assumption of an additional symmetry of the polytope (Thm. \ref{symmy}). 
Furthermore we show that in fixed dimension 
a centrally symmetric simplicial reflexive polytope with the maximal number of vertices 
is even uniquely determined (Thm. \ref{csymmy}). For the proofs we generalise a result due to Casagrande 
\cite[Thm. 2.4]{Cas03a} to the case of a $\Q$-factorial Gorenstein toric Fano variety (Cor. \ref{pic-bound}) 
saying that the Picard number of a 
nonsingular toric Fano variety exceeds the Picard number of a torus-invariant prime divisor at most by three. 

In section 6 we give sharp bounds on the number of lattice points in terminal Fano polytopes (Cor. \ref{termpts}) 
and centrally symmetric reflexive polytopes (Thm. \ref{central}). 

\begin{acknowledgement}{\rm
The author would like to thank his thesis advisor 
Professor Victor Batyrev for posing problems, his advice 
and encouragement, as well as Professor G\"unter Ewald for giving 
reference to \cite{Wir97} and Professor Klaus Altmann for the possibility of giving a talk at the FU Berlin. 
The author would also like to thank Professor Maximillian Kreuzer for 
the support with the computer package PALP, the classification data and many examples.

The author was supported by DFG, Forschungsschwerpunkt "Globale Methoden in der komplexen Geometrie". 
This work is part of the author's thesis.}
\end{acknowledgement}

\section{Notation and basic definitions}
\label{sec:1}

In this section we fix the notation for toric varieties and lay out the basic notions of toric Fano varieties. 
We refer to \cite{Ewa96} for combinatorial convexity and \cite{Ful93,Oda88} for toric varieties. 
In \cite{Bat94} reflexive polytopes were introduced. For a survey of (toric) Fano 
varieties see \cite{Deb01}. 

Let $N \cong \Z^d$ be a $d$-dimensional lattice and $M = \Hom_\Z(N,\Z) \cong \Z^d$ the dual lattice with 
$\pro{\cdot}{\cdot}$ the nondegenerate symmetric pairing. 
As usual, $\NQ = N \otimes_\Z \Q \cong \Q^d$ and $\MQ = M\otimes_\Z \Q \cong \Q^d$ 
(respectively $\NR$ and $\MR$) will denote the rational (respectively real) scalar extensions. 

{\em Throughout the paper the roles of $N$ and $M$ are interchangeable.}

For a subset $S$ in a real vector space let $\lin(S)$ (respectively $\aff(S)$, $\conv(S)$, $\pos(S)$) be the 
linear (respectively affine, convex, positive) hull of $S$. 
A subset $P \subseteq \MR$ is called a polytope, if it is the convex hull of finitely many points in $\MR$. 
The boundary of $P$ is denoted by $\randp$, the relative interior of $P$ by $\relint P$. When 
$P$ is full-dimensional, its relative interior is also denoted by $\intr P$. A face $F$ of $P$ is denoted by 
$F \leq P$, the vertices of $P$ form the set $\V(P)$, the facets of $P$ the set $\F(P)$. 
$P$ is called a lattice polytope, respectively rational polytope, 
if $\V(P) \subseteq M$, respectively $\V(P) \subseteq \MQ$. An isomorphism of 
lattice polytopes is an isomorphism of the lattices such that the induced real linear isomorphism maps 
the polytopes onto each other. 

We usually denote by $\fan$ a complete fan in $\MR$. 
The $k$-dimensional cones of $\fan$ form a set $\fan(k)$. The elements in $\fan(1)$ are called rays, 
and given $\tau \in \rays$, we let $\vt$ denote the unique generator of $M \cap \tau$.

{\em Since the polytope $P$ is throughout contained in $\MR$, fans will also have cones mostly in $\MR$ 
in contrast to the usual convention.} 

There are two possibilities to define a complete fan from a polytope $P \subseteq \MR$ that is rational 
and $d$-dimensional:

First we can define the normal fan $\calN_P$ of $P$, i.e., an element of $\calN_P$ 
is the closed cone of inner normals of a face of $P$.

Second if the origin is contained in the interior of $P$, i.e., $0 \in \intr P$, then we can 
define $\Sigma_P := \{\pos(F) \,:\, F \leq P\}$ as the fan that is spanned by $P$. 
In this case there is the important notion of the dual polytope
\[P^* := \{y \in \NR \,:\,\, \pro{x}{y} \geq -1 \;\forall\, x \in P\},\]
that is a rational $d$-dimensional polytope with $0 \in \intr P^*$. We have $\Sigma_{P^*} = \calN_P$. 

Duality means $(P^*)^* = P$. 
There is a natural combinatorial correspondence between $i$-dimensional faces of $P$ and 
$(d-1-i)$-dimensional faces of $P^*$ that reverses inclusion. 
For a facet $F \leq P$ we let $\eta_F \in \NQ$ denote the uniquely determined 
inner normal satisfying $\pro{\eta_F}{F} = -1$, so 
$\V(P^*) = \{\eta_F \,:\, F \in \F(P)\}$.

The dual of the product of $d_i$-dimensional polytopes $P_i \subseteq \R^{d_i}$ with $0 \in \intr P_i$ for $i=1,2$ is 
given by
\begin{equation}(P_1 \times P_2)^* = \conv(P_1^* \times \{0\}, \{0\} \times P_2^*) \subseteq \R^{d_1} \times \R^{d_2}.
\label{product}
\end{equation}

\begin{definition}{\rm
Let $X$ be a normal complex variety and $\KX$ the {\em canonical divisor} of $X$, i.e., 
a Weil divisor of $X$ whose restriction to the regular locus defines the canonical sheaf there.

$\KX$ is {\em $\Q$-Cartier}, if there exists a positive integer $j$ such that $j \KX$ is a Cartier divisor. 
The smallest such $j$ is called the {\em Gorenstein index} $\jX$ of $X$. 

A complex variety $X$ is called {\em Fano variety} (respectively {\em weak 
Fano variety}), if $X$ is projective, normal and 
the anticanonical divisor $-\KX$ is an ample (respectively nef and big) $\Q$-Cartier divisor.

A complex variety $X$ is called {\em Gorenstein}, iff $\jX = 1$, i.e., $\KX$ is a Cartier divisor.}
\end{definition}

Back to the toric case: Let $\fan$ be a complete fan in $\MR$, we denote by $X := X(M,\fan)$ the associated 
complete normal toric variety. 

$X$ is projective if and only if there is a 
$d$-dimensional lattice polytope $P \subseteq \MR$ with $0 \in \intr P$ and $\fan = \Sigma_P$, or equivalently, 
if there exists a $d$-dimensional lattice polytope $Q \subseteq \NR$ with $\fan = \calN_Q$.

Set $\Gf := \{\vt \::\: \tau \in \rays\}$ and $G_\sigma := \Gf \cap \sigma$ for $\sigma \in \fan$.
 
It is well-known that $\KX := - \sum_{\tau \in \rays} \Vt$ is a canonical divisor of $X$, 
where $\Vt$ denotes the torus-invariant prime divisor associated to the ray $\tau$. 
Furthermore $\KX$ is $\Q$-Cartier if and only if for all $\sigma \in \fan$ the set $\Gf \cap \sigma$ 
is contained in an affine hyperplane. 

Now we define the lattice polytope $\Pf := \conv(\Gf) \subseteq \MR$ with 
$0 \in \intr \Pf$, and the rational dual polytope $\Qf := P_\fan^* \subseteq \NR$. Then 
$X$ is a {\em toric Fano variety} if and only if $\fan = \Sigma_\Pf$, or equivalently, $\fan = \calN_{\Qf}$. 
The Gorenstein index $\jX$ is the minimal $k$ such that $k \Qf$ is a lattice polytope. 

Whenever $\KX$ is $\Q$-Cartier we can explicitly compute the {\em discrepancy} of $X$ 
(see \cite[Prop. 12]{Deb01}):
\[\disc(X) = -1 + \min\{ \pro{\us}{v} \::\: \sigma \in \fan(d),\, 0 \not= v \in \sigma \cap (M \backslash \Gs)\},\]
thus a rational number in $]-1,1]$. {\em Especially $X$ has log-terminal singularities}, i.e., $\disc(X) > -1$.

The following definitions are now convenient: 

\begin{definition}{\rm
Let $P \subseteq \MR$ be a $d$-dimensional lattice polytope with $0 \in \intr P$.
\begin{itemize}
\item $P$ is called a {\em Fano polytope}, if the vertices are primitive lattice points.
\item $P$ is called a {\em canonical Fano polytope}, if $\intr P \cap M = \{0\}$.
\item $P$ is called a {\em terminal Fano polytope}, if $P \cap M = \{0\} \cup \V(P)$.
\item $P$ is called a {\em smooth Fano polytope}, if the vertices of any facet of $P$ form a $\Z$-basis of 
the lattice $M$.
\end{itemize}}
\end{definition}

\begin{remark}{\rm
{\em Beware:} In most papers (see \cite{Bat99,Sat00}) a Fano polytope is assumed to be already 
a smooth Fano polytope. This more systematic notation was partly introduced in \cite{Deb01}.}
\end{remark}

We have the following correspondence theorem:

\begin{proposition}
There is a correspondence between isomorphism classes of Fano polytopes and isomorphism classes of 
toric Fano varieties.

Thereby canonical Fano polytopes correspond to toric Fano varieties with canonical singularities, i.e., $\disc X \geq 0$; 

terminal Fano polytopes correspond to toric Fano varieties with terminal singularities, i.e., $\disc X > 0$;

smooth Fano polytopes correspond to nonsingular toric Fano varieties.
\end{proposition}

There is the following important finiteness theorem (see \cite{Bor00} for a survey):

\begin{theorem}
For $\epsilon > 0$ there exist only finitely many isomorphism classes of $d$-dimensional toric Fano varieties 
with discrepancy greater than $-1 + \epsilon$.
\label{finiteness}
\end{theorem}

For the weak case we define:

\begin{definition}{\rm
Let $P \subseteq \MR$ be a Fano polytope spanning $\fan$. Then a fan $\fans$ is called a {\em crepant refinement} 
of $\fan$, if $\fans$ is a refinement of $\fan$ in the usual sense and additionally 
for any $\taus \in \rayss$ there exists a $\sigma \in \fan$ such that $\vts \subseteq \conv(\Gs)$. 
When the toric variety associated to the fan $\fans$ is again projective, the fan $\fans$ is called a 
{\em coherent crepant refinement}.}
\end{definition}

Using the ramification formula (see \cite{Deb01}) we see that such crepant refinements correspond to 
equivariant proper birational morphisms\newline$f : X\str = X(M,\fans) \to X = X(M,\fan)$ with $K_{X\str} = f^* \KX$. 

\begin{proposition}
Toric weak Fano varieties correspond uniquely up to isomorphism to coherent crepant refinements of fans spanned 
by Fano polytopes.
\end{proposition}

Finally let's recall the following definition:

\begin{definition}{\rm
A complex normal variety is called {\em $\Q$-factorial}, if any Weil divisor is $\Q$-Cartier.}
\end{definition}

For the next result (see \cite[Thm. 2.2.24]{Bat94}) recall that a polytope is simplicial, 
if any facet is a simplex:

\begin{proposition}
$\Q$-factorial toric Fano varieties 
correspond uniquely up to isomorphism to simplicial Fano polytopes.

There exists a coherent crepant refinement by stellar subdivisions 
that resolves a weak toric Fano variety $X$ with canonical singularities to a $\Q$-factorial weak toric Fano variety 
$X\str$ with terminal singularities. 
\label{reso}
\end{proposition}

Such a morphism $X\str \to X$ is called a {\em MPCP-desingularization} in \cite{Bat94}. 
$X$ is said to {\em admit a coherent crepant resolution}, if such an $X\str$ can be chosen to be nonsingular.

Reflexive polytopes naturally enter the picture when considering Gorenstein toric Fano varieties.

\begin{definition}{\rm
A Fano polytope whose dual is a lattice polytope is called {\em reflexive polytope}.}
\end{definition}

The following proposition is now straightforward.

\begin{proposition}
Gorenstein toric Fano varieties correspond uniquely up to isomorphism to reflexive polytopes.

Gorenstein toric weak Fano varieties correspond uniquely up to isomorphism to 
coherent crepant refinements of fans spanned by reflexive polytopes.
\end{proposition}

There is a complete {\em duality} of reflexive polytopes and their duals.

\newpage
\begin{proposition}
Let $P \subseteq \MR$ be a $d$-dimensional lattice polytope with $0 \in \intr P$. 
Then the following conditions are equivalent:
\begin{enumerate}
\item $P$ is a reflexive polytope
\item $P$ is a lattice polytope and $P^*$ is a lattice polytope
\item $P^*$ is a reflexive polytope
\end{enumerate}
If this holds, then $\intr P \cap M = \{0\}$, i.e., $P$ is a canonical Fano polytope.
\end{proposition}

The proof is straightforward by using the following local property that characterises reflexive polytopes:

\begin{lemma}
Let $P \subseteq \MR$ be a reflexive polytope.

For any $F \in \F(P)$ and $m \in F \cap M$ there is a $\Z$-basis $e_1, \ldots, e_{d-1}, e_d := m$ of $M$ such that 
$F \subseteq \{x \in \MR \,:\, x_d = 1\}$; thus $\eta_F = - e_d^*$ in the dual $\Z$-basis $e_1^*, \ldots, e_d^*$ of $N$. 

\label{root-local}
\end{lemma}

\begin{proof}

Since $\eta_F \in N$, the following short exact sequence splits
$$0 \to \{x \in M \,:\, \pro{\eta_F}{x} = 0\} \to M \stackrel{\pro{\eta_F}{\cdot}}{\to} \Z \to 0.$$

\end{proof}

As an immediate corollary from \ref{finiteness} we get:

\begin{corollary}
Gorenstein toric (weak) Fano varieties have canonical singularities. 
In fixed dimension $d$ there is only a finite number of isomorphism classes of $d$-dimensional 
Gorenstein toric Fano varieties.
\end{corollary}

Finally refering to \ref{reso} there is following result (see \cite{Bat94} and \ref{bsp}):

\begin{proposition}
For $d\leq 3$ any $d$-dimensional Gorenstein toric Fano variety admits a coherent crepant resolution.
\label{crep3}
\end{proposition}

For the proof we recall a convenient definition:

\begin{definition}{\rm
A lattice polytope $P \subseteq \MR$ is called {\em empty}, if $P \cap M = \V(P)$.}
\end{definition}

Now \ref{crep3} follows from \ref{reso} and the next well-known lemma:

\begin{lemma}
Let $P$ be a $d$-dimensional lattice polytope with $0 \in \intr P$.
\begin{enumerate}
\item Let $d = 2$. 
Lattice points $x,y$ form a $\Z$-basis if and only if $\conv(0,x,y)$ is an empty two-dimensional polytope. 
If $P$ is a canonical Fano polytope, this is fulfilled, if $x,y$ are lattice points on the boundary that are 
not contained in a common facet of $P$ and $x + y \not= 0$.

In particular a two-dimensional terminal Fano polytope is a smooth Fano polytope, and 
a two-dimensional canonical Fano polytope is a reflexive polytope.
\item Let $d = 3$ and $P$ be reflexive. 
Three lattice points $x,y,z$ in a common facet of $P$ form a $\Z$-basis 
if and only if $\conv(x,y,z)$ is empty.

In particular a three-dimensional simplicial terminal reflexive polytope is a smooth Fano polytope.
\end{enumerate}
\label{klein}
\end{lemma}

\section{Projecting along lattice points on the boundary}
\label{sec:2}

{\em In this section $P$ is a $d$-dimensional reflexive polytope in $\MR$.}

The projection map along a vertex of $P$ is an essential tool in investigating 
toric Fano varieties, since 
one hopes to get some information from the corresponding lower-dimensional variety (see \cite{Bat99,Cas03a}). 
In the case of a reflexive polytope it is also useful to consider projecting along lattice points on the boundary 
of $P$ that are not necessarily vertices.

The following definitions will be used throughout the paper:

\begin{definition}{\rm
Let $x,y \in \randp$ with $x \not= y$.

\begin{itemize}
\item $[x,y] := \conv(x,y)$,\; $]x, y] := [x,y] \backslash \{x\}$,\; $]x,y[ := [x,y] \backslash \{x,y\}$. 
\item $x \sim y$, if $[x,y]$ is contained in a face of $P$, i.e., 
$x$ and $y$ are contained in a common facet of $P$. 
\item The {\em star set} of $x$ is the set 
\[\st(x) := \{y \in \randp \,:\, x \sim y\} = \bigcup \{F \in \F(P) \,:\, x \in F\}.\]

\item The {\em link} of $x$ is the set
\[\randx := \rand\, \st(x) = \bigcup \{G \leq P \,:\, G \subseteq \st(x),\, x \not\in G\}.\]

\item $y \vdash x$ ($y$ is said to be {\em away from $x$}), if $y$ is not in the relative interior of $\st(x)$. 
Hence
\[\randx = \{y \in \st(x) \,:\, y \vdash x\}.\]
$y$ is away from $x$ iff there exists a facet that contains $y$ but not $x$, e.g., if $y$ is a vertex or 
$x \not\sim y$. There is also a local criterion:
$$ y \vdash x \;\;\lolra\;\; x + \lambda (y-x) \not\in P \,\;\forall\: \lambda > 1.$$ 

\end{itemize}
\label{away}}
\end{definition}

The next proposition shows the important properties of the projection of a reflexive polytope along some lattice point 
on the boundary.

\begin{proposition}
Let $P \subseteq \MR$ be a reflexive polytope. 
Let $v \in \randPM$. Define the quotient lattice $M_v := M / \Z v$ and the canonical projection map along v
\[\Pi = \Pi_v \,:\, M_\R \to (M_v)_\R = M_\R / \R v.\]
Then $P_v := \Pi_v(P)$ is a lattice polytope in $(M_v)_\R$ with $\V(P_v) \subseteq \Pi(\V(P))$ containing the origin 
in the interior. 

\begin{enumerate}
\item Let $U$ be the set of elements $x \in P$ such that 
$x + \lambda v \not\in P$ for all $\lambda > 0$. Then we have a canonical bijection 
\[U \to P_v.\] 
We denote the inverse map by $\rho$. 

For $S := \st(v)$ we have 
\[U = S\]
and thus
\[P_v = \conv(\Pi(\V(P) \cap \rand v)).\]

\item The projection map induces a bijection $S \cap M \to P_v \cap M$.

\item The projection map induces a bijection $\rand v \to \rand P_v$. 
$\rand P_v$ is covered by the projection of all $(d-2)$-dimensional faces $C$ such that 
$C \in \F(F)$ for some facet $F$ of $P$ with $v \in F$ and $v \not\in C$; $\Pi(C)$ is contained in a facet of $P_v$.

\item\[\rho(\V(P_v)) \subseteq \V(P) \cap \rand v.\]

Let $z \in \V(P) \cap \rand v$. Then $\Pi(z) \in \rand P_v$. 
If $]v,z[$ is contained in the relative interior of a facet of $P$, then $\Pi(z) \in \V(P_v)$.

\item Let $F \in \F(P)$ with $v \in F$. Then
\[\pos(\Pi(F)) \cap P_v = \Pi(F).\]

\item The image $\Pi(F)$ of a facet $F$ parallel to $v$, i.e., $\pro{\eta_F}{v} = 0$, is a facet of $P_v$. It 
is $\Pi^{-1}(\Pi(F)) \cap P = F$. There are at least $\card{\V(\Pi(F))}$-vertices of $F$ in $S$. 
Any point in $F \cap S$ is contained in a facet that contains $v$ and intersects $F$ in a $(d-2)$-dimensional face.

The preimage $\Gamma := \Pi^{-1}(F') \cap P$ of a facet $F'$ of $P_v$ is either a facet of $P$ parallel to $v$ 
or a $(d-2)$-dimensional face of $P$. In the last case $\Gamma \to F'$ is an isomorphism, and 
there exists exactly one facet of $P$ that contains $\Gamma$ and $v$.

\item Suppose $-v \in P$. Then any facet of $P$ either contains $v$, or $-v$, or is parallel to $v$, 
i.e., a facet of the form $\Pi^{-1}(F') \cap P$ for $F' \in \F(P_v)$. 

\item \[(P_v)^* \cong P^* \cap v^\perp \text{ as lattice polytopes}.\]
$P_v$ is reflexive if and only if $P^* \cap v^\perp$ is a lattice polytope.

\end{enumerate}

\label{proj}
\end{proposition}

\begin{proof}

1. Let $F$ be a facet of $P$ containing $v$ and $x$. If $\lambda > 0$, 
then $\pro{\eta_F}{x + \lambda v} = -1 - \lambda < -1$, so $x + \lambda v \not\in P$. Hence $S \subseteq U$. 

On the other hand let $x \in U$. Considering the polytope $P \cap \lin(v,x)$ we see there is 
a facet $F$ of $P$ not parallel to $v$ that contains $x$ with $\pro{\eta_F}{v} < 0$. 
Since $P$ is reflexive, we have $\pro{\eta_F}{v} = -1$, hence $v \in F$ and $x \sim v$. 
This implies $S=U$.

2. Let $m' \in P_v \cap M_v$. We have $u := \rho(m') \in U=S$. So there exists a facet $F \in \F(P)$ 
with $v,u \in F$. By \ref{root-local} there is a $\Z$-basis $e_1, \ldots, e_{d-1}, e_d=v$ of $M$ such that 
$\V(F) \subseteq \{x \in \MR \,:\, x_d = 1\}$, i.e., $\eta_F$ is the dual vector $-e^*_d$. Let $u = 
\lambda_1 e_1 + \ldots + \lambda_d e_d$ for $\lambda_1, \ldots, \lambda_d \in \R$. Now $e_1, \ldots, e_{d-1}$ 
is a $\Z$-basis of $M_v$. Therefore $\lambda_1, \ldots, \lambda_{d-1} \in \Z$. Since $u \in F$, we get 
$\lambda_d = u_d = \pro{-\eta_F}{u} = 1$, hence $u \in M$. 

3. Let $z \in \rand v$. Then $z \vdash v$ and there exists a facet $F \in \F(P)$ with $v,z \in F$. 
Assume $z' := \Pi(z) \not\in \rand P_v$. Then there is some $\lambda > 1$ such that 
$\lambda z' \in \rand P_v$. We have $x := \rho(\lambda z') \sim v$. Since $x \in \lin(v,z)$ and $z \sim v$, this 
yields $x \in F$. This implies that $z$ is a proper convex combination of $v$ and $x$, a contradiction 
to $z \vdash v$. The remaining properties are as easily seen.

4. The first statements follow from the first and the third point. 
Now suppose $z \in \V(P) \cap S$ such that there is only one facet $F \in \F(P)$ that contains $v$ and $z$. Since 
$z \in \V(F)$ we can choose an affine hyperplane $H$ that 
intersects $F$ only in $z$ and is parallel to $v$. For $P_v := \Pi_v(P)$ and $H' := \Pi_v(H)$ let $x' \in H' \cap P_v$. 
It remains to show that $x' = z' := \Pi_v(z)$. So assume not. $H'$ intersects $\Pi_v(F)$ only in $z'$. 
Therefore $\rho(y') \not\in F$ for all $y' \in ]z',x']$. Finiteness of $\F(P)$ implies that $z$ is contained 
in another facet $\not= F$ containing $v$, a contradiction. 

5. This is proven as the third point.

6. The first statements follow from the third and the fourth point.

For the second statement 
let $\dim(\Gamma) = d-2$. Now observe that if $\Gamma \to F'$ were not injective, a facet containing $\Gamma$ 
necessarily would be parallel to $v$, so its image a facet containing $F'$, a contradiction. Therefore 
$\Gamma \to F'$ is an isomorphism of polytopes with respect to their affine hulls. 
Now choose $x \in \relint F'$. Let $y := \rho(x) \in S \cap \Gamma$. By assumption also $y \in \relint \Gamma$. 
Let $G \in \F(P)$ with $v,y \in G$. Then $\Gamma \subseteq G$, hence $G$ is one of the two facets containing $\Gamma$.

7. Let now $-v \in P$. Any facet $F \in \F(P)$ satiesfies $-1 \leq \pro{\eta_F}{v} = - \pro{\eta_F}{-v} \leq 1$. 
From this the statements follow.

8. Choose again a facet $F$ of $P$ with $v \in F$, and a $\Z$-basis 
$e_1, \ldots, e_d$ of $M$ such that $e_d = v$ and 
$e_1, \ldots, e_{d-1}$ is a $\Z$-basis of $M \cap \eta_F^\perp \cong M_v$. In these coordinates we get 
$(P_v)^* = \{y' \in \R^{d-1} \,:\, \pro{y'}{x'} \geq -1 \; \forall x' \in P_v\} = 
\{y' \in \R^{d-1} \,:\, \pro{(y'_1, \ldots, y'_{d-1}, 0)}{(x_1, \ldots, x_d)} \geq -1 \; \forall x \in P\} 
= P^* \cap v^\perp$.

The last point and its proof are taken from \cite[Thm. in Sect. 3]{AKMS97}.

\end{proof}

Let's consider the {\em algebraic-geometric interpretation} of the projection map:

Let $\fan := \Sigma_P$, $X := X(M,\Sigma_P)$. Let $v \in \V(P)$, $\tau := \pos(v) \in \rays$ and $\Pi = \Pi_v$. 
As in \cite[1.7]{Oda88} the fan $\fanqu := \{\Pi(\sigma) \::\: \sigma \in \fan, \tau \leq \sigma\} = 
\{\pos(\Pi(F)) \::\: F \in \F(P), v \in F\}$ defines the projective toric variety $\Vt$ that is the 
torus invariant prime divisor corresponding to the ray $\tau$. 

On the other hand there is the projected polytope $P_v := \Pi(P)$ that 
spans a fan $\fanv$ in $(M_v)_\R$, we denote the corresponding projective toric variety by $X_v$. 
In the following we discuss how and when $\Vt$ and $X_v$ are related.

We choose a triangulation $\T := \{T_k\}$ of $\rand v$ into simplicial lattice polytopes. 
Then $\fanT := \{\pos(\Pi(T_k))\}$ is a simplicial fan in $(M_v)_\R$ with corresponding $\Q$-factorial 
complete toric variety $\XT$. From prop. \ref{proj}(5) it follows that $\fanT$ is a {\em common refinement} 
of $\fanqu$ and $\fanv$. Especially there are induced proper birational morphisms 
$\XT \to \Vt$ and $\XT \to X_v$.

In general $\fanv$ is not a refinement of $\fanqu$. However in the case that $P$ is simplicial, we can 
choose $\T$ obviously in such a way that $\fanT = \fanqu$, in particular this implies that 
$\fanqu$ is a refinement of $\fanv$. 

In order to draw conclusions about the canonical divisor and singularities of these lower-dimensional 
toric varieties there is the following sufficient assumption:
\begin{equation}
\exists\: f \in \N_{>0} \;\;:\;\; \forall\; w \in \VP \cap \rand v \;:\; \card{[v,w] \,\cap\, M} - 1 = f.
\label{fcond}
\end{equation}
Suppose this condition holds. For any $w \in \VP \cap \rand v$ the second point in the proposition implies 
$\card{[0,\Pi_v(w)] \cap M_v} -1 = f$. Furthermore let $w' \in \V(P_v)$, $\tau' := \pos(w')$. 
Since prop. 
\ref{proj}(4) implies $\rho(w') \in \V(P) \cap \rand v$, the previous consideration yields $\vtg = (1/f) w'$, 
hence $P_{\fanv} = (1/f) P_v$ (as defined in section \ref{sec:1}). Therefore 
$-K_\Xv$ is an ample $\Q$-Cartier divisor, i.e., $\Xv$ is a {\em toric Fano variety}. 

\begin{definition}{\rm
$X$, respectively $P$, is called {\em semi-terminal}, if for all $v \in \V(P)$ condition 
(\ref{fcond}) holds for $f=1$, i.e., $[v,w] \cap M = \{v,w\}$ for all $w \in \VP \cap \rand v$.
\label{semi-terminal}}
\end{definition}

\begin{proposition}
Let $P \subseteq \MR$ be a reflexive polytope. 
\begin{enumerate}
\item 
$P$ is semi-terminal iff $P_v$ is a Fano polytope for all $v \in \VP$.
\item $P$ is terminal iff $P_v$ is a canonical Fano polytope for all $v \in \VP$.
\end{enumerate}
\label{termprop}
\end{proposition}

\begin{proof}

1. From left to right: This holds since $P_\fanv = P_v$ is a Fano polytope.

From right to left: Let $v \not= w \in \V(P)$ with $v \sim w$. Choose $C$ and $F$ as in prop. \ref{proj}(3) such 
that $w \in \V(C)$. By prop. \ref{proj}(3) we see that 
$F' := \aff(\Pi_v(C)) \cap P_v$ is a facet of $P_v$. Hence by prop. \ref{proj}(6) 
for $\Gamma := \Pi_v^{-1}(F') \cap P$ we have either $\Gamma = G$ for a facet 
$G \in \F(P)$ that is parallel to $v$ or $\Gamma = C$. In the first case $\pro{\eta_G}{v} = 0$ and 
$\pro{\eta_G}{w} = -1$, so $\card{[v,w] \cap M} = 2$. 
In the second case obviously $F' = \Pi_v(C)$, so 
$\Pi_v(w) \in \V(F')$, hence by assumption a primitive lattice point. From prop. \ref{proj}(2) we get 
$\card{[v,w] \cap M} = 2$. 

2. From left to right: Let $v \in \VP$, $0 \not= m' \in M_v \cap P_v$. In the notation of prop. \ref{proj}(1,2) 
we get $\rho(m') \in M \cap \randp \subseteq \VP$ by assumption. 
Hence $m' = \Pi_v(\rho(m')) \in \rand P_v$ by prop. \ref{proj}(4). 

From right to left: Assume there is a $w \in \randPM$, $w \not\in \VP$. Then $w$ is a proper convex 
combination of vertices of $P$ contained in a common facet. 
Let $v$ be one of them. Then $\Pi_v(w)$ is obviously in the interior of $P_v$, a 
contradiction.

\end{proof}

\begin{corollary}
Let $P \subseteq \MR$ be a reflexive polytope.\\Then the following two statements are equivalent:
\begin{enumerate}
\item $X$ has terminal singularities
\item $X$ is semi-terminal and $\Xv$ has canonical singularities for any $v \in \V(P)$
\end{enumerate}
If this holds, then $\Xv$ is a toric Fano variety for any $v \in \V(P)$.
\label{corXv}
\end{corollary}

In particular we see that terminality is a necessary condition for obtaining a reflexive polytope under projection, 
however not sufficient for $d \geq 4$. 

\begin{example}{\rm
Let $d=4$, and $e_1, \ldots, e_4$ a $\Z$-basis of $M$. We define the simplicial centrally symmetric reflexive polytope 
$P := \conv(\pm (2e_1 + e_2 + e_3 + e_4), \pm e_2, \pm e_3, \pm e_4)$. Then $P$ is combinatorially 
a crosspolytope, has $8$ vertices and $16$ facets. It is a terminal Fano polytope 
but not a smooth Fano polytope, so especially it admits no 
crepant resolution. The projection $P_{e_4}$ along the vertex $e_4$ has $6$ vertices, $P_{e_4}$ is even a terminal 
Fano polytope but not reflexive. 
This polytope is taken from \cite{Wir97} where it is used in a different context.
\label{bsp}}
\end{example}

Now let's look at $\Vt$: To ensure that the canonical divisor of $\Vt$ is $\Q$-Cartier, we need 
in general the $\Q$-factoriality of $X$. So let $P$ be simplicial and assume again that condition (\ref{fcond}) holds. 
Then $\fanqu$ is a coherent crepant refinement of $\fanv$. 
Hence $-K_\Vt$ is a nef $\Q$-Cartier divisor, i.e., $\Vt$ is a {\em toric weak Fano variety}. 
From cor. \ref{corXv} we get:

\pagebreak
\begin{corollary}
Let $X$ be a $\Q$-factorial Gorenstein toric Fano variety.\\The following two statements are equivalent:
\begin{enumerate}
\item $X$ has terminal singularities
\item $X$ is semi-terminal and $\Vt$ has terminal singularities for any $\tau \in \rays$
\end{enumerate}
If this holds, then $\Vt$ is a $\Q$-factorial toric weak Fano variety for any $\tau \in \rays$.
\label{corVt}
\end{corollary}

Finally to additionally derive the Gorenstein property, i.e., that the canonical divisor is $\Z$-Cartier, 
we need a stronger assumption, that is trivial in the case of a smooth Fano polytope:

\vspace{-\bigskipamount}

\begin{gather}
\text{For any $F \in \F(P)$ with $v \in F$ and $C \in \F(F)$ with $v \not\in C$}\notag\\
\text{there exist $w_1, \ldots, w_{d-1} \in C \cap M$ such that}\label{scond}\\
\text{$w_1 - v, \ldots, w_{d-1} - v$ is a $\Z$-basis (of $\eta_F^\perp \cap M$).}\notag
\end{gather}

If this condition is fulfilled, then (\ref{fcond}) holds for $f=1$, 
$P_v$ is reflexive by prop. \ref{proj}(3), and $\Xv$ is a {\em Gorenstein toric Fano variety}. 
If $P$ is also simplicial, then $\Vt$ is a {\em Gorenstein toric weak Fano variety}. 

Suppose now $X$ is semi-terminal, simplicial and $\Vt$ is nonsingular for any $\tau \in \rays$. 
It follows from \ref{root-local} that (\ref{scond}) holds for any $v \in \V(P)$. 
Then $P_v$ is reflexive, in particular canonical for any $v \in \VP$, hence cor. \ref{corXv} 
implies that $P$ is terminal. Since $P$ is also simplicial, the assumption implies that 
$P$ is already a smooth Fano polytope. We have proven the following corollary:

\begin{corollary}
Let $X$ be a $\Q$-factorial Gorenstein toric Fano variety.\\The following two statements are equivalent:
\begin{enumerate}
\item $X$ is nonsingular
\item $X$ is semi-terminal and $\Vt$ is nonsingular for any $\tau \in \rays$
\end{enumerate}
If this holds, then $\Vt$ is a $\Q$-factorial toric weak Fano variety for any $\tau \in \rays$, and 
$\Xv$ is a Gorenstein toric Fano variety admitting the coherent crepant resolution $\Vt \to \Xv$ 
for any $v \in \VP$, $\tau = \pos(v)$.
\label{nonsi}
\end{corollary}

The important fact that the projection of a smooth Fano polytope is reflexive was 
already proven by Batyrev in \cite[Prop. 2.4.4]{Bat99}, however he used 
the notion of a primitive relation \cite{Bat91} and results of Reid about the Mori cone \cite{Rei83}. 

There is now a generalisation of this result to the class of toric Fano varieties with locally complete 
intersections. These varieties were thoroughly investigated in \cite{DHZ01}, where it was proven that 
they admit coherent crepant resolutions.

\begin{proposition}
Let $X$ be a Gorenstein toric Fano variety that has singularities that are locally complete intersections.\\
Then the following three statements are equivalent:
\begin{enumerate}
\item $X$ is semi-terminal
\item $X$ has terminal singularities
\item Any facet of $P$ can be embedded as a lattice polytope in $[0,1]^{d-1}$
\end{enumerate}

If this holds, then $\Xv$ is a Gorenstein toric Fano variety for any $v \in \VP$. 
If additionally $X$ is $\Q$-factorial, then $X$ is nonsingular.
\end{proposition}

\begin{proof}

The facets of the corresponding reflexive polytope $P$ are so called {\em Nakajima} polytopes, 
a comprehensive description can be found in \cite{DHZ01}. Using their results it is 
straightforward to prove the following statement by induction on $n$:

Let $v$ be a vertex of an $n$-dimensional Nakajima polytope $F$ in a lattice $M$ such that 
$\card{[w,w'] \cap M} = 2$ for all $w,w' \in \V(P)$, $w \not= w'$. 
Then $F$ is empty, can be embedded as a lattice polytope in $[0,1]^n$, and 
for any facet $C \in \F(F)$ with $v \not\in C$ there exist $n$ vertices $w_1, \ldots, w_n$ of $C$ such that 
$w_1 - v, \ldots, w_n - v$ is a $\Z$-basis.

From this result the proposition is obvious using condition (\ref{scond}).

\end{proof}

\section{Pairs of lattice points on the boundary}
\label{sec:3}

{\em Throughout the section let $P$ be a $d$-dimensional reflexive polytope in $\MR$.}

In \cite{Bat91} the important notions of {\em primitive collections} and 
{\em primitive relations} were defined for nonsingular projective 
toric varieties, and used in \cite{Bat99} for the classification of four-dimensional smooth Fano polytopes. 
Unfortunately these useful tools cannot simply be generalised to the class of reflexive polytopes. 
However the next proposition shows that in the simplest yet most important case of a primitive collection 
of order two, i.e., a pair of lattice points on the boundary that is not contained in a common face, 
we still have a kind of generalised primitive relation:

\begin{proposition}
Let $P \subseteq \MR$ be a reflexive polytope, $v, w \in \randPM$, $v \not= w$. 

Exactly one of the following three statements holds:

\begin{enumerate} 
\item $v \sim w$
\item $v + w = 0$
\item $v + w \in \randp$
\end{enumerate}

Let the third condition be fulfilled. Then it holds:

$v,w$ is a $\Z$-basis of $\lin(v,w) \cap M$. There exists exactly one pair $(a,b) \in \N_{>0}^2$ with 
$z := z(v,w) := a v + b w \in \randp$ such that $v \sim z$ and $w \sim z$. Moreover:
\begin{enumerate}
\item[{\rm i.}] $a=1$ or $b=1$. $a =\: \card{[w,z] \cap M} - 1$ and $b =\: \card{[v,z] \cap M} - 1$.\\
If $F \in \F(P)$ with $v,z \in F$, then $\pro{\eta_F}{w} = \frac{a-1}{b}$.
\item[{\rm ii.}] Any facet containing $z$ (or $v+w$) contains exactly one of the points $v$ or $w$.
\item[{\rm iii.}] For any $F \in \F(P)$ containing $v$ and $z$ there exists a facet $G \in \F(P)$ containing 
$w$ and $z$ such that $F \cap G$ is a $(d-2)$-dimensional face of $P$.
\item[{\rm iv.}] If $z \in \V(P)$, $b=1$, and $]v,z[$ is contained in 
the relative interior of a facet of $P$, then $[w,z]$ is contained in an edge.
\end{enumerate}
\label{prim}
\end{proposition}

\begin{proof}

Let $v \not\sim w$ and $v + w \not= 0$. The first condition implies that for any facet $F \in \F(P)$ we have 
$\pro{\eta_F}{v+w} = \pro{\eta_F}{v} + \pro{\eta_F}{w} > -2$. However reflexivity of $P$ 
implies that this must be a natural number greater or equal to $-1$, so $v+w \in P$ by duality. We get 
$0 \not= v+w \in \randp$, because $P$ is canonical. 

$v,w$ is a $\Z$-basis by \ref{klein}(1).

Let $F$ be a facet of $P$ containing $v+w$. We may assume 
$\pro{\eta_F}{v} = -1$ and $\pro{\eta_F}{w} = 0$. This implies $v \sim v+w$. 
We can use this consideration again for the pair $v+w,w$. 
Since $F \cap M$ is finite, this eventually yields a natural number 
$b \in \N_{>0}$ such that $z = v + b w \in F$ and $w \sim z$. In particular $a=1$. This proves 
the existence of $z$ and ${\rm i}$.

${\rm ii}$. Let $F' \in \F(P)$ with $z \in F'$. 
Assume $v,w \not\in F'$, hence $-1 = \pro{\eta_{F'}}{z} = a \pro{\eta_{F'}}{v} + b \pro{\eta_{F'}}{w} \geq 0$, 
a contradiction. 

${\rm iii}$.
Let $F \in \F(P)$ with $v,z \in F$. We set $w' := (a-1)v+w \in \randp \cap M$, 
hence $z=z(v,w')$. Since $\pro{\eta_F}{w'} = 0$ and $z \in F \cap \st(w')$, 
prop. \ref{proj}(6) applied to $\Pi_{w'}$ implies that there exists a facet $G \in \F(P)$ containing $z$ and $w'$ 
that intersects $F$ in a $(d-2)$-dimensional face. Since obviously $v \not\in G$, it follows from the 
second point that $w' \in G$. This immediately yields $w \in G$.

${\rm iv}$. Follows from prop. \ref{proj}(4) applied to $\Pi_v$.

\end{proof}

The next figure illustrates the proposition:

\begin{center}\epsfxsize1.25in\epsffile{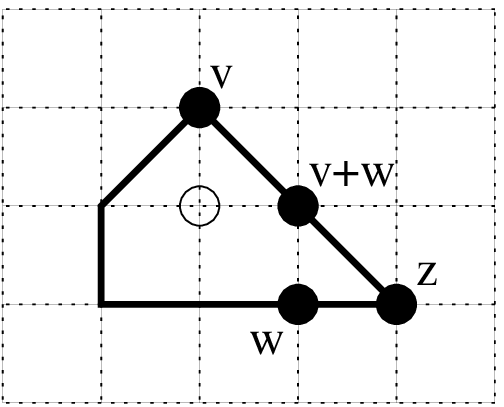}\end{center}

The symmetric relation $\sim$ defines a {\em graph} $\W(P)$ on $\randPM$. 
From the previous proposition we can now easily derive the following corollary about combinatorial 
properties of reflexive polytopes:

\begin{corollary}
\label{comb}
Let $P \subseteq \MR$ be a reflexive polytope.
\begin{enumerate}
\item Any pair of points in $\randPM$ can be connected by at most three edges 
of the graph $\W(P)$, with equality only possibly occuring for a centrally symmetric pair of points. 
In particular the diameter of the graph $\W(P)$ is at most three.
\item The previous statements also hold for the restriction of $\W(P)$ to the set of vertices, 
which is a purely combinatorial object. 
In the case of a simplicial polytope this is just the usual edge-graph on the vertices of $P$.
\item By dualizing we get that a pair of facets of a reflexive polytope is either parallel, 
contains a common vertex, or does have mutually non-trivial intersection with another facet. 
\end{enumerate}
\end{corollary}

Without using the existing classification of two-dimensional reflexive polytopes (see prop. \ref{two}) 
the proposition and the corollary yield an immediate application in the case of $d=2$ (for the proof of the second point 
use statement {\rm i} of the proposition).

\begin{corollary}
Let $P$ be a two-dimensional reflexive polytope.
\begin{enumerate}
\item $P$ has at most six vertices; equality occurs iff $P$ is of type $6a$ in prop. \ref{two}.
\item Any facet of $P$ contains at most five lattice points; there exists a facet with five lattice points 
iff $P$ is of type $8c$ in prop. \ref{two}.
\end{enumerate}
\label{six}
\end{corollary}

This first point is also a direct consequence of \cite[Thm. 1]{PR00} saying that $\card{\rand P \cap M} + 
\card{\rand P^* \cap N} = 12$, where however no direct combinatorial proof is known that does not use some kind 
of induction. 

Another application is to show that certain combinatorial isomorphismtypes 
of polytopes cannot be realized as reflexive polytopes. As an example have a look at the {\em regular polyhedra} 
(see for instance \cite{Sti01}). 

\begin{corollary}
In any dimension the combinatorial type of a $d$-simplex, a $d$-cube, and a $d$-crosspolytope 
can be realized as a reflexive polytope. 
There is no three-dimensional reflexive polytope that is combinatorially isomorphic to a 
dodecahedron or an icosahedron. There is no four-dimensional reflexive polytope that is 
combinatorially isomorphic to the 120-cell or 600-cell.
\label{noiso}
\end{corollary}

\begin{proof}

The first statement is trivial.

Let $P$ be a reflexive polytope and $d=3$. By duality we can assume that $P$ is 
combinatorially isomorphic to an icosahedron. Cor. \ref{comb}(2) yields that $P$ is centrally symmetric. 
However any three-dimensional 
centrally symmetric simplicial reflexive polytope has at most $8$ vertices as will be proven in thm. \ref{csymmy}.

Finally by cor. \ref{comb} and duality it is enough to note that the diameter of 
the edge-graph of the simplicial 600-cell is larger than three (see \cite[Fig. 5]{Sti01}).

\end{proof}

It is now an astonishing observation (see \cite{KS00b}) 
that the self-dual 24-cell {\em can} be uniquely realized as a reflexive polytope 
with vertices $\{\pm e_i \,:\, i = 1, \ldots, 4\}$\\$\cup \{\pm (e_i - e_j) \,:\, 
i = 1,2,\, j = 2,3,4,\, j > i\} \cup \{\pm (e_i - e_3 - e_4) \,:\, i = 1,2\} \cup 
\{\pm (e_1 + e_2 - e_3 - e_4)\}$ for $e_1, \ldots, e_4$ a $\Z$-basis of $M$. 
It is even centrally symmetric and terminal. Here it is interesting to note the necessity of these conditions:

\begin{corollary}
Let $P$ be a four-dimensional reflexive polytope $P$ that is combinatorially a 24-cell. Then $P$ 
has to be centrally symmetric and terminal.
\end{corollary}

\begin{proof}

Let $v$ be a vertex of $P$. Now choose the vertex $w \in \V(P)$ corresponding to the usual antipodal point. 
Assume $v + w \not= 0$. It is easy to see (see \cite[Fig. 4]{Sti01}) 
that the intersection of a facet containing $v$ and a facet containing $w$ is empty or consists of a unique vertex $z$ 
where $]v,z[$ and $]w,z[$ are contained in the relative interiors of these facets. 
This implies $v \not\sim w$ and $z(v,w) = z \in \VP$, a contradiction to the last point of prop. \ref{prim}.

The terminality of $P$ can be proven in an analogous way.

\end{proof}

\pagebreak
\section{Classification results in low dimensions}
\label{sec:4}

Smooth Fano polytopes, as they form the most important class of reflexive polytopes, 
were intensively studied over the last decade by Batyrev 
\cite{Bat91,Bat99}, Casagrande \cite{Cas03b,Cas03a}, Debarre \cite{Deb01}, 
Sato \cite{Sat00}, et al. It could be rigorously proven that there are 
$18$ smooth Fano polytopes for $d=3$ (see \cite{Bat82,WW82}) and $124$ for $d=4$ (see \cite{Bat99,Sat00}) up to 
isomorphism. 

Here we will have a look at recent classification results of reflexive polytopes in low dimensions.

For $d=1$ the polytope $[-1,1]$ corresponding to $\P^1$ is the only Fano polytope. 
For $d=2$ any canonical Fano polytope is reflexive by \ref{klein}(1), 
and these isomorphism classes can be easily classified (e.g., see \cite{KS97} or \cite[Thm. 6.22]{Sat00}). 
For the convenience of the reader and later reference we will give the list of the $16$ isomorphism classes of reflexive 
polygons as well as a simple proof.

\begin{proposition}
There are exactly $16$ isomorphism classes of two-dimen\-sional reflexive polytopes 
(the number in the labels are the number of lattice points on the boundary):

\begin{center}\epsfxsize3.75in\epsffile{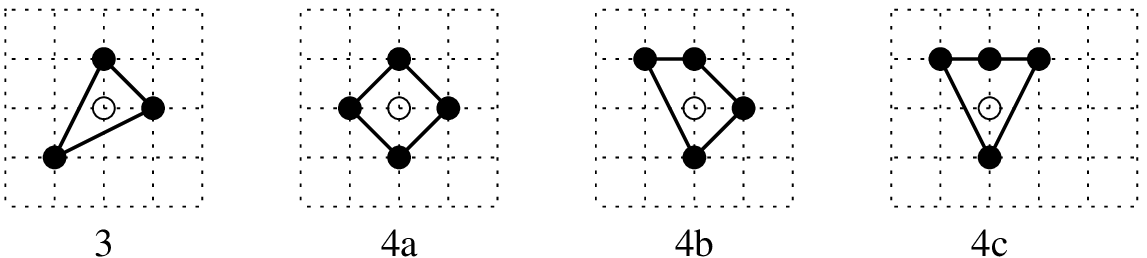}\end{center}
\begin{center}\epsfxsize3.75in\epsffile{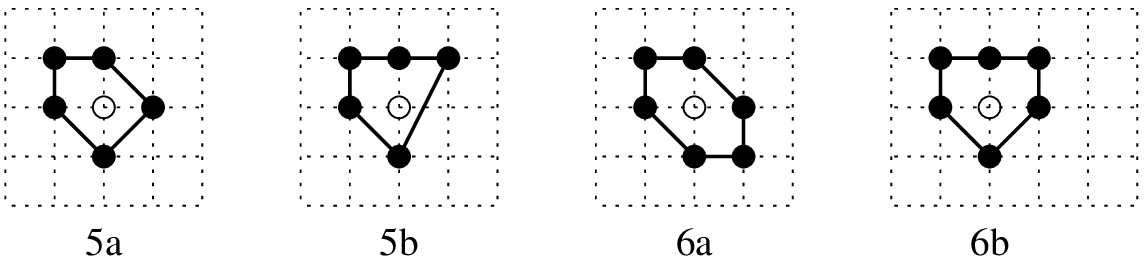}\end{center}
\begin{center}\epsfxsize3.75in\epsffile{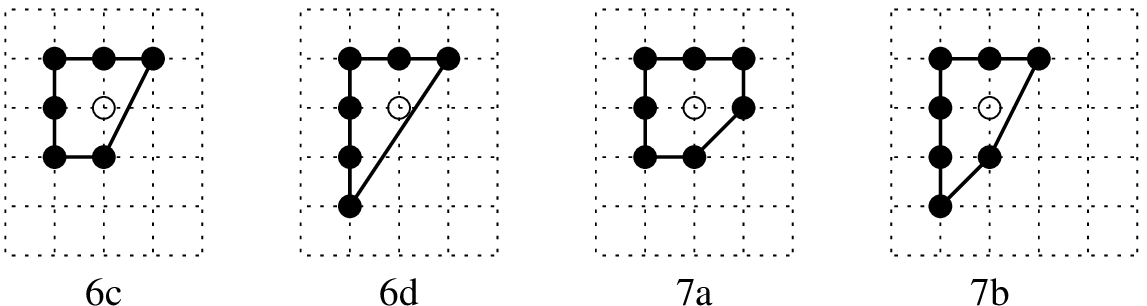}\end{center}
\begin{center}\epsfxsize3.75in\epsffile{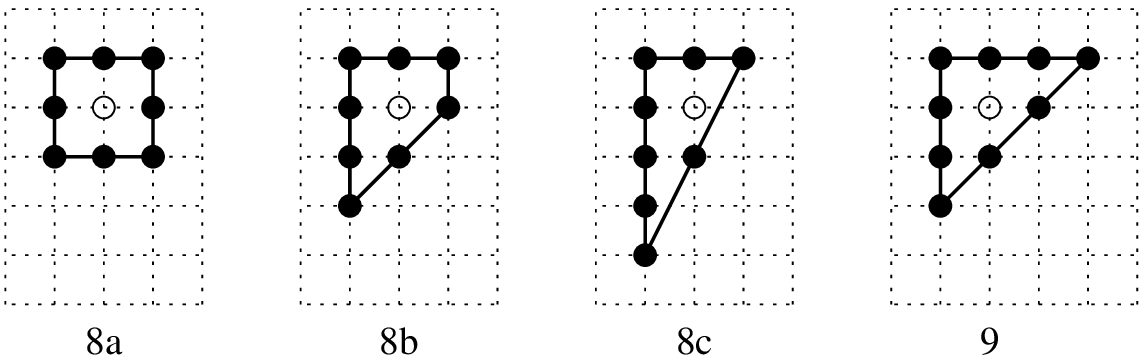}\end{center}

\label{two}
\end{proposition}

\begin{proof}

Let $P$ be a two-dimensional reflexive polytope. We distinguish three different cases:

\begin{enumerate}
\item Any facet of $P$ contains only two lattice points, i.e., $P$ is a terminal Fano polytope. 
There are three different cases (see prop. \ref{prim}):

\begin{enumerate}
\item $P$ is combinatorially a triangle. 

By \ref{klein}(1) we may assume that $(1,0)$, $(0,1)$ are vertices of $P$. 
Let $x$ be the third vertex. From prop. \ref{termprop}(2) it follows that the projection of $P$ along $(1,0)$ 
is a canonical Fano polytope, 
i.e., isomorphic to $[-1,1]$, hence $x_2 = -1$. By projecting along $(0,1)$ we get $x_1 = -1$, 
so $P$ is of type $3$.

\item There exist three vertices $u,v,w \in \VP$ with $u+w = v$. 

Since $P$ is a terminal Fano polytope, prop. \ref{prim} implies that $u \sim v$ and $w \sim v$, and 
we may assume $u = (-1,1)$, $v = (0,1)$, $w = (1,0)$. Again projecting along $v$ yields 
$P \cap \{(-1,x) \,:\, x \in \Z\} \subseteq \{(-1,0), (-1,1)\}$,
$P \cap \{(0,x) \,:\, x \in \Z\} \subseteq \{(0,-1), (0,0), (0,1)\}$,
$P \cap \{(1,x) \,:\, x \in \Z\} \subseteq \{(1,-1), (1,0)\}$. 
We get as possible types $4b$, $5a$, $6a$.

\item Any two vertices that are no neighbours are centrally symmetric. 

This immediately implies that $P$ is of type $4a$.
\end{enumerate}

\item There exists a facet $F$ containing exactly one lattice point in $\relint F$. 

We may assume $\V(F) = \{(-1,1), (1,1)\}$. 
Then by prop. \ref{proj}(1) we have $P \subseteq \{x \in \R^2 \,:\, -1 \leq x_1 \leq 1,\, x_2 \leq 1\}$. 
Since $(0,-1)$ is not contained in $\intr P$, we get $P \subseteq \{x \in \R^2 \,:\, x_2 \geq -3\}$. 
From this we readily derive the next ten isomorphism types $4c$,$5b$,$6b$,$6c$,$6d$,$7a$,$7b$,$8a$,$8b$,$8c$.

\item The remaining case.

We may assume $\V(F) = \{(-1,1), (a,1)\}$ for $a \in \N$, $a \geq 2$. Let $v \in \VP$ with $v_2 \leq -1$ minimal. 
As $(1,0)$ is not in the interior of $P$, we have $v_1 \leq 0$. 
Then by assumption necessarily $\conv((-1,1), (-1,-2), (2,1)) \subseteq P$. This must be an equality, 
hence $P$ is of type $9$.
\end{enumerate}
\end{proof}

The proof includes the statement that there are exactly five toric Del Pezzo surfaces, i.e., two-dimensional nonsingular 
toric Fano varieties, a result that can also be proven by birational factorisation \cite[Prop. 2.21]{Oda88}, 
primitive relations \cite{Bat91} or determinants \cite[Thm. V.8.2]{Ewa96}.

In general even for $d=3$ there are too many reflexive polytopes to give a classification 
by pencil and paper. In \cite{Con02} reflexive simplices are classified using the notion of weights. 
Kreuzer and Skarke described in \cite{KS97,KS98,KS00a} a general algorithm to classify reflexive polytopes 
in fixed dimension $d$. Using their computer program PALP (see \cite{KS02}) they applied their method for $d \leq 4$, 
and found $4319$ reflexive polytopes for $d=3$ and $473800776$ for $d=4$. They also described how to find a 
normal form of lattice polytopes, toric fibrations and symmetries.

However it is still interesting to find rigorous mathematical proofs of observations and classification results of 
smaller classes of reflexive polytopes by directly using their intrinsic properties. 
There is the following result by the author:

\begin{theorem}
There are $100$ three-dimensional terminal reflexive polytopes. 
A three-dimensional reflexive polytope is terminal 
if and only if any facet is either a simplex 
where the vertices form a $\Z$-basis or a parallelogram $\conv(v,w,x,y)$ where 
$v,x,w$ form a $\Z$-basis and satisfy the equation $v+x=w+y$.
\label{100}
\end{theorem}

A complete proof is contained in the thesis of the author. It relies on the notion of an {\em AS-point}, i.e., 
a vertex that is both additive, i.e., the sum of two other vertices, and symmetric, 
i.e., its antipodal point is also a vertex. 
If no such AS-point exists, we use prop. \ref{termprop}(2) to show that the polytope has at most eight vertices, so by 
\ref{crep3} we can use the classification 
of three-dimensional proper nonsingular toric varieties with Picard number five or less 
which are minimal in the sense of equivariant blow-ups as described in \cite[1.34]{Oda88}. On the other hand if 
there exists an AS-point, we can use the tools in the previous sections (in particular lemma \ref{klein}(2), 
prop. \ref{proj}(7) 
and prop. \ref{prim}) to completely describe the polytope by a suitably generalised notion of a primitive relation.

Recently Kasprzyk classified in \cite{Kas03} all $634$ three-dimensional terminal polytopes by first describing 
the minimal cases purely mathematically and then using a computer program for the remaining ones.

\section{Sharp bounds on the number of vertices}
\label{sec:5}

{\em Throughout the section let $P$ be a $d$-dimensional reflexive polytope in $\MR$.}

In higher dimensions only in very special cases classification results exist. So one tries to 
find at least sharp bounds on invariants and to characterise the case of equality. Here we examine the 
number of vertices of a reflexive polytope. 

This invariant corresponds to the rank of the class group of the associated toric variety $X := X(M,\Sigma_P)$. 
There is the exact sequence (see \cite[Prop. 3.4]{Ful93})
\[0 \to M \to \DivTX \to \Cl(X) \to 0,\]
where $\DivTX$ denotes the free abelian group of torus invariant Weil divisors. Since $\DivTX$ has rank 
$\card{\Sigma_P(1)} = \card{\V(P)}$, we have
\[\rank \,\Cl(X) = \card{\V(P)} - d.\]

The classification of Kreuzer and Skarke shows that the maximal number of vertices of a $d$-dimensional reflexive 
polytope is $6$ for $d=2$, $14$ for $d=3$ and $36$ for $d=4$. This observation motivates a conjecture:

\begin{definition}{\rm
We define $\Zt_2 := \conv(\pm [0,1]^2)$, it is the (up to isomorphism) unique centrally-symmetric 
self-dual smooth Fano polytope with $6$ vertices (of type $6a$ in prop. \ref{two}). 
We denote by $S_3 := X(M,\Sigma_{\Zt_2})$ the associated nonsingular toric Del Pezzo surface that 
is the blow up of $\P^2$ in three torus-invariant points.}
\end{definition}

\begin{conjecture}
Let $P$ be a $d$-dimensional reflexive polytope. Then
\[\card{\V(P)} \leq 6^{\frac{d}{2}},\]
where equality occurs if and only if $d$ is even and $P \cong (\Zt_2)^{\frac{d}{2}}$. 
\label{vertconj}
\end{conjecture}

\begin{remark}{\rm
{\em It would be enough to prove this conjecture for $d$ even}, because products of reflexive polytopes are again 
reflexive.}
\end{remark}

This conjecture is confirmed by the computer classification for $d \leq 4$, however 
only for $d=2$ a rigorous proof is known (see cor. \ref{six}(1) and thm. \ref{bounds}(1) below).

The next theorem yields two coarse upper bounds on the number of vertices of a reflexive polytope in terms 
of some combinatorial invariants of the facets. The first bound is a straightforward generalisation 
of a bound due to Voskresenskij and Klyachko \cite[Thm.1]{VK85} originally proven 
in the setting of a smooth Fano polytope. The second upper 
bound is a generalisation of \cite[Thm. 8]{Deb01}, where Debarre improved from a bound of 
order $O(d^2)$ on 
the number of vertices of a smooth Fano polytope to a bound of order $O(d^{3/2})$, which is the 
asymptotically best upper bound that is known at the moment. 
We recover the original results for simplicial reflexive polytopes.

\begin{theorem}
Let $P$ be a reflexive polytope. 

Define $\alpha := \max( \V(F) \,:\, F \in \F(P))$ and $\beta := \max( \F(F) \,:\, F \in \F(P))$.

\begin{enumerate}
\item $\card{\V(P)} \leq 2 d \alpha$.

More precisely we distinguish two cases:

If $\alpha \geq 2 d - 3$, then $\card{\V(P)} \leq 2 d (\alpha - d + 2) - 2$.

If $\alpha \leq 2 d - 3$, then $\card{\V(P)} \leq d \alpha + \alpha - d + 1$.

If $P$ is simplicial, i.e., $\alpha = d$, and $d \geq 3$, this yields
\[\card{\V(P)} \leq d^2 + 1.\]

\item $\card{\V(P)} \leq (\alpha - d + 1) \beta + 2 + 2 \sqrt{(\alpha-1) (d+1) ((\alpha-1) + (\alpha-d+1) \beta)})$.

If $P$ is simplicial, i.e., $\alpha = d = \beta$, this yields
\[\card{\V(P)} \leq d + 2 + 2 \sqrt{(d^2-1)(2d-1)}.\]
\end{enumerate}
\label{bounds}
\end{theorem}

\begin{proof}

Analysing the proofs of Thm. 1 in \cite{VK85} and Thm. 8 in \cite{Deb01} in the more general 
setting of a reflexive polytope, we see that by taking the general invariants $\alpha$ and $\beta$ into account we just 
have to reprove remark 5(2) in section 2.3 of \cite{Deb01}, because only there explicitly a lattice basis was used. 
That result is essentially the first part of the next lemma.

\end{proof}

\begin{lemma}
Let $P \subseteq \MR$ be a reflexive polytope. 

Let $F \in \F(P)$, $u := \eta_F \in \V(P^*)$ and $\{F_i\}_{i \in I}$ the facets that 
intersect $F$ in a $(d-2)$-dimensional face. Let $m \in \randp \cap M$ with $\pro{u}{m} = 0$. 

Then $m \in \cup_{i \in I} F_i$. 

Let additionally $F$ be a simplex with 
$\V(F) = \{e_1, \ldots, e_d\}$. Let $e_1^*, \ldots, e_d^*$ be the dual $\R$-basis of $\NR$. 
For $i = 1, \ldots, d$ denote by $F_i$ the facet of $P$ 
such that $F_i \cap F = \conv(e_j \,:\, j \not=i)$ and choose a lattice point $m^i$ on $F_i$ 
that is not contained in $F$.
\begin{enumerate}
\item For $i \in \{1, \ldots, d\}$ we have
\[m \not\in F_i \;\lolra\; \pro{e_i^*}{m} \geq 0.\]
\item If there exists $i \in \{1, \ldots, d\}$ such that $m \in F_i$ and 
$m \not\in F_j$ for all $j \in \{1, \ldots, d\}$, $j\not=i$, then $m \not\sim e_i$.
\item If $\pro{u}{m^i} = 0$ and $\pro{e_i^*}{m^i} = -1$ for $i = 1, \ldots, d-1$, 
then $e_1, \ldots, e_d$ is a $\Z$-basis of $M$.
\end{enumerate}
\label{fund}
\end{lemma}

\begin{proof}

The first part follows from prop. \ref{proj}(6) for $\Pi_m$.

Now let $F$ be a simplex. Then $u = \sum_{j=1}^d (-e_j^*) \in N$. Let $i \in \{1, \ldots, d\}$. 
Since $m^i  \not\in F$ 
and $0$ is in the interior of $P$, the number $\alpha_i := \frac{-1-\pro{u}{m^i}}{\pro{e_i^*}{m^i}} > 0$ is well-defined. 
We get $\eta_{F_i} = u + \alpha_i e_i^*$. From this 1. is readily derived. 2. is just a corollary. 
In 3. we get $\alpha_i = 1$ and $e_i^* = \eta_{F_i} - u \in N$ for 
$i = 1, \ldots, d-1$ and $e_d^* = -\eta_{F_i} - e_1^* - \ldots - e_{d-1}^* \in N$.

This proof is inspired by remark 5(2) in section 2.3 of \cite{Deb01}.

\end{proof}

In the following we will focus on the class of {\em simplicial} reflexive polytopes, i.e., 
where the associated varieties are $\Q$-factorial, or equivalently, the class number equals the Picard number. 
The previous theorem already gave a hint that simplicial reflexive polytopes 
are actually quite close to smooth Fano polytopes at least when considering only the number of vertices. 
Also in this case there is an explicit conjecture:

\begin{conjecture}
Let $P$ be a $d$-dimensional simplicial reflexive polytope. Then 
\[\card{\V(P)} \leq \left\{\begin{array}{lcl}3 d &,& d \text{ even,}\\3 d - 1 &,& d \text{ odd.}\end{array}
\right.\]
For $d$ even equality holds if and only if $P^* \cong (\Zt_2)^{\frac{d}{2}}$, i.e., $X \cong (S_3)^{\frac{d}{2}}$.
\label{mainconj}
\end{conjecture}

\begin{remark}{\rm
{\em It would be enough to prove this conjecture for $d$ even:} Assume there were a simplicial reflexive 
polytope $P$ with $d$ odd and $\card{\V(P)} \geq 3d$. Then necessarily $P^* \times P^* \cong (\Zt_2)^d$, 
this would imply $P$ to be centrally symmetric with $\card{\V(P)} = 3d$, a contradiction to $d$ odd.}
\label{enough}
\end{remark}

The bound is also sharp in the odd-dimensional case, take $X = \P^1 \times (S_3)^{\frac{d-1}{2}}$. 
However even for $d=3$ there is exactly one another simplicial reflexive polytope with $8$ vertices, it is a 
smooth Fano polytope, not centrally symmetric, and the associated toric variety $X$ 
is an equivariant $S_3$-fibre bundle over $\P^1$. 

From \cite{Kas03} we get that reflexivity is essential, because 
the maximal number of vertices a three-dimensional simplicial terminal Fano polytope can have 
is $10$.

This conjecture was originally proposed in the case of smooth Fano polytopes by Batyrev, 
and was rigorously proven to hold for (up to) five-dimensional smooth Fano polytopes 
by Casagrande in \cite[Thm. 3.2]{Cas03a}. 
In the above form the conjecture is confirmed by the computer classification of Kreuzer and Skarke for $d\leq 4$. 
Moreover 
it yields that there are $194$, respectively $5450$, classes of three-dimen\-sional, respectively four-dimensional, 
simplical reflexive polytopes; however only $151$ four-dimensional terminal simplicial reflexive polytopes. 

The main goal of this section is to give a proof in the case of additional symmetries of the polytope.

\begin{theorem}
Conjecture \ref{mainconj} holds in the case of a simplicial reflexive polytope $P$ where the dual polytope $P^*$ 
contains a vertex $u \in \V(P^*)$ such that $-u \in P^*$.
\label{symmy}
\end{theorem}

\begin{theorem}
Let $P$ be a centrally symmetric simplicial reflexive polytope.

If $d$ is even, then $\card{\V(P)} \leq 3 d$, with equality iff $P^* \cong (\Zt_2)^{\frac{d}{2}}$.

If $d$ is odd, then $\card{\V(P)} \leq 3 d - 1$, with equality iff $P^* \cong [-1,1] \times (\Zt_2)^{\frac{d-1}{2}}$.

\label{csymmy}
\end{theorem}

We need some preparation for the proofs.

The main result for analysing smooth Fano polytopes is a theorem of Reid about 
extremal rays of the Mori cone and primitive relations (see \cite{Rei83} and \cite[Thm. 1.3]{Cas03a}). 
Although for simplicial reflexive polytopes there is no general notion of a primitive relation, for the simplest 
case as defined in prop. \ref{prim} we still have an analogous result (recall definition \ref{away}):

\begin{lemma}
Let $P$ be a simplicial reflexive polytope. 

Let $v \in \V(P)$, $w \in \randp \cap M$ with $v+w \in \randp$ and $z := z(v,w)$. 

Let $x \in \randp$, $x \not\in \{v,w,z\}$, with $x \sim z$ and $x$ away from $v$. 

Then $\conv(x,z,w)$ is contained in a face.

Moreover exactly one of the following two conditions holds:

\begin{enumerate}
\item Any facet containing $x$ and $z$ contains also $w$.
\item There exists a facet $F$ with $x,v,z \in F$. 
\end{enumerate}

The second case must occur, if $w \in \V(P)$ and $x$ is away from $w$.

If the second case occurs, we have:

For any such $F$ there exists a unique facet $G$ with $x,w,z \in G$ such that 
$F \cap G$ is a $(d-2)$-dimensional face of $P$. 
$F\cap G$ consists of those elements of $F$ that are away from $v$, respectively those elements of $G$ that 
are away from the (unique) vertex not in $F$. Obviously $w \not\in F$ and $v \not\in G$. 

\label{reid}
\end{lemma}

\begin{proof}

Assume the first case is not fulfilled. Prop. \ref{prim}{({\rm ii})} implies that there exists 
a facet $F \in \F(P)$ with $x,z,v \in F$. By \ref{prim}{({\rm iii})} there is a facet $G \in \F(P)$ 
containing $w$ and $z$ such that $F \cap G$ is a $(d-2)$-dimensional face. Since $F,G$ are simplices and 
$v \in \V(F)$, the remaining statements are now straightforward.

\end{proof}

The next result is a generalisation of a lemma proven by Casagrande \cite[Lemma 2.3]{Cas03a} 
for smooth Fano polytopes, 
here we recover the original statement in the more general setting of a terminal simplicial reflexive polytope.

\begin{lemma}
Let $P$ be a simplicial reflexive polytope. 

Let $v,w \in \V(P)$, $w' \in \randPM$ away from $w$. Furthmore let $v+w \in \randp$ and $v+w' \in \randp$, 
$z := z(v,w)$ and $z' := z(v,w')$. 

We define $K := P \cap \lin(v,w,w')$. 
Then $K$ is a {\em two-dimensional} reflexive polytope (of possible types $5a$, $6a$, $6b$, $7a$ in prop. \ref{two}).

If $K$ is terminal, then $z=v+w$, $w'=-v-w=-z$, $z'=v+w'=-w$; and {\em either} 
$\rand K \cap M = \{v,w,z,w',z'\}$ {\em or} $\rand K \cap M = \{v,w,z,w',z',-v = w + w' = z(w,w')\}$.

\label{casa}
\end{lemma}

\begin{proof}

Let $z = a v + b w$ and $z' = a' v + b' w'$ as in \ref{prim}. 
We note that $w'$ and $z$ is away from $v$ and $w$; $z'$ is away from $v$ and $w'$.

Assume $w' \sim z$. Since 
$w'$ is away from $w \in \V(P)$, it follows from \ref{reid} that there exists a facet that 
contains $w',z,v$, hence $w' \sim v$, a contradiction.

Thus $w' \not\sim z$.

There are now two different cases, and it must be shown that the second one cannot occur.
\begin{enumerate}
\item $v,w,w'$ are linearly dependent.

By \ref{prim} there are three possibilities:

If $w \sim w'$, then $K = \conv(v,z,w,z',w')$. 
If $w + w' = 0$, then $v \in \conv(v+w,v+w')$, a contradiction. 
If $v+w \in \randp$, then $K = \conv(v,z,w,z',w',$\\$z(w,w'))$.

Thus in any case $K$ is a lattice polytope with $5$ or $6$ vertices, 
canonical, hence by \ref{klein}(1) reflexive. By analysing 
the cases in \ref{two} we get the remaining statements.

\item $v,w,w'$ are linearly independent.

Hence also $z,z',v$ are linearly independent. 

Assume $z' \sim z$. 
\ref{reid} implies that $\conv(z',z,w)$ is contained in a facet $F \in \F(P)$. Since $v \not\in F$ and $z' \in F$, 
\ref{prim}({\rm ii}) implies $w' \in F$, a contradiction to $w' \not\sim z$.

Thus $z \not\sim z'$.

By assumption $z+z'\not=0$, hence $z+z' \in \randp$. Let $y := z(z',z) = k z' + l z \in \randPM$. 
We have $y \not\in \{z,z',v,w,w'\}$, because $v,w,w'$ are linearly independent. 
Choose $y' \in \randp$ with $y' = v + \lambda (y - v)$ for $\lambda \geq 1$ maximal, so that 
$y'$ is away from $v$. Furthermore $y' \sim z'$ and $y' \sim z$. So by \ref{reid} there exist 
facets $F_1, F_2 \in \F(P)$ such that $\conv(y',z,w) \subseteq F_1$ and $\conv(y',z',w') \subseteq F_2$; 
$v \not\in F_1, F_2$. 

Now choose $y'' = w + \mu (y' - w) \in P$ for $\mu \geq 1$ maximal; so $y''$ is away from $w$. Furthermore 
$\conv(y'',y',z,w) \subseteq F_1$ and 
$y''$ away from $v$, so by \ref{reid} there exists a facet $G \in \F(P)$ that contains $y'',v,z$ and intersects 
$F_1$ in a $(d-2)$-dimensional face. Hence necessarily 
$\pro{\eta_{F_1}}{v} = \frac{b-1}{a}$ and $\pro{\eta_G}{w} = \frac{a-1}{b}$.

$K$ is a three-dimensional polytope. Any face of $K$ is contained in 
a face of $P$. We have $F_1 \not= F_2$, since $w' \not\sim z$. 
So $C := F_1 \cap F_2 \cap K$ is a vertex or edge of $K$ containing $y'$. 
Since also $w' + z \not= 0$ and $w' \not\sim z$, we get $w' + z \in \randp$. 
We set $x := z(w',z) \in \randPM$. 

We distinguish several cases:

\begin{enumerate}
\item $y' = y$ .

\begin{enumerate}
\item $w \in F_2$.

The vertices of $C$ consist of $y''$ and $w$. It is also easily seen that $x \in C$ with $w \not= x \not=y$, 
so also $x \not=y''$. If $a=1$, it were 
$0 = \pro{\eta_G}{w} > \pro{\eta_G}{x} > \pro{\eta_G}{y''} = -1$, a contradiction. 

Therefore $a \geq 2$, $b=1$, $\pro{\eta_G}{w} = a-1 > 0$.

Then $\pro{\eta_{F_2}}{v} = \frac{b'-1}{a'}$ and 
$-1 = \pro{\eta_{F_2}}{y} = -k + l a \frac{b'-1}{a'} - l b$, thus, since $b=1$,
\begin{equation}\frac{k-1}{l} = a \frac{b'-1}{a'} - 1 \in \N.\label{eq2}\end{equation}

Since $\frac{b'-1}{a'} \not= 0$, $b' \geq 2$ and $a' = 1$. 

If $k=1$, then (\ref{eq2}) yields $1 = a (b'-1) \geq a \geq 2$, a contradiction.

If $l=1$, then $\pro{\eta_G}{y} = - k -1 + k b' \pro{\eta_G}{w'}$. So 
$-1 = \pro{\eta_G}{y''} = (1-\mu) (a-1) + \mu \pro{\eta_G}{y} = (1-\mu) (a-1) 
- \mu k - \mu + \mu k b' \pro{\eta_G}{w'}$. This implies 
$\pro{\eta_G}{w'} = \frac{\frac{-a}{\mu}+k+a}{kb'}$. 
Since (\ref{eq2}) yields $k = a (b'-1)$, this implies 
\[\pro{\eta_G}{w'} = \frac{\frac{-1}{\mu b'}+ 1}{b'-1} \in \N.\]

On the other hand $\mu \geq 1$ and $b' \geq 2$ yields $0 < \frac{-1}{\mu b'}+ 1 < 1$, this contradicts the 
previous equation.

\item $w \not\in F_2$.

This immediately implies $y'' = y$. We find $x' \in F_1$ such that $x \in [w,x']$ and 
$x'$ is away from $w$. By \ref{reid} we have $x' \in F_1 \cap G$. As $w' \not\sim z$, 
we obtain $x'=y$. Hence $x \in ]w,y[$. 

Since $w' \sim x$ there exists a facet $H \in \F(P)$ containing $w',w,y$; furthermore $H \not=F_2$, 
since $w \not\in F_2$. 
Hence there are edges $F_2 \cap H \supseteq [w',y]$ and $F_1 \cap H = [w,y]$ of $K$. 

Since $w \not\sim z'$ we can define in a double recursion $x_r^0 := x$, $x_l^0 := x(z',w)$, $x_r^i := x(x_l^{i-1},z)$, 
$x_l^i := x(z', x_r^{i-1})$ for $i \in \N$, $i \geq 1$. As $w',y,z$ and $w,y,z'$ are linearly independent, 
we easily see that this procedure is well-defined, and 
$x_l^0, x_l^1, x_l^2, \ldots$ are pairwise different lattice points in $]w',y[$ and 
$x_r^0, x_r^1, x_r^2, \ldots$ are pairwise different lattice points in $]w,y[$. Hence we have 
constructed infinitely many lattice points in $P$, a contradiction.

\end{enumerate}

\item $y' \not= y$.

If $y'' \not= y'$, then obviously $y$ is a lattice point in the interior of $K$, a contradiction.

Thus $y'' = y'$. This implies $y \in ]v,y'[$, so $y \in G$, hence 
$\conv(z',y',y,v)$ is contained in a facet $F' \in \F(P)$. \ref{reid} implies that there exists 
a unique facet $G' \in \F(P)$ that contains $w',z',y'$ such that $F' \cap G'$ is a $(d-2)$-dimensional face. 

Furthermore $\frac{b-1}{a} = \pro{\eta_{F_1}}{v} > \pro{\eta_{F_1}}{y} > \pro{\eta_{F_1}}{y'} = -1$, hence 
$b\geq 2$ and $a=1$. Especially we get $\pro{\eta_G}{w} = 0$.

This yields again $x \in F_1$.

\begin{enumerate}
\item $w \in G'$.

Since $w \in \V(P)$, $w \not\in F'$ and $w'$ is away from $w \in \V(G')$, \ref{reid} implies that 
$w' \in G' \cap F'$, a contradiction.

\item $w \not\in G'$.

Let $C' := G' \cap F_1 \cap K$. Then $y'$ is a vertex of $C'$.

Assume $C'$ were an edge. Let $v' \in \V(C')$ with 
$v' \not= y'$. This implies $v' \not= w$. Then $v'$ is away from $w$, hence by \ref{reid} $v' \in G$, therefore 
$v' = y'$, a contradiction.

So $C' = \{y'\}$, and just the same way we see that $[y',w]$ is an edge of $K$. Therefore 
$x \in ]w,y']$. Since $0 = \pro{\eta_G}{w} > \pro{\eta_G}{x} \geq \pro{\eta_G}{y'} = -1$, this implies $x = y'$.

Furthermore we have 
$y' = (1-\lambda) v + \lambda y = ((1-\lambda) + \lambda (k a' + l a)) v + \lambda k b' w' + \lambda l b w$, 
where $\lambda > 1$. On the other hand $x = r w' + s z = r w' + s a v + s b w$ for $r,s \in \N$, $r,s \geq 1$. 

Comparing the coefficients for $w'$ and $w$ this yields
$$\lambda l = s, \;\;\; \lambda k b' = r.$$
From the first equation we get $s = \lambda l > l \geq 1$, so $s \geq 2$. This implies $1 = r = \lambda k b' 
> k b' \geq 1$, a contradiction.
\end{enumerate}
\end{enumerate}
\end{enumerate}
\end{proof}

Using prop. \ref{proj}(1-4) and analysing the possible cases 
in prop. \ref{two} it is straightforward to prove a corollary of the previous lemma:

\begin{corollary}
Let $P$ be a simplicial reflexive polytope and $v \in \V(P)$. 

There are at most three vertices of $P$ not in the star set of $v$; equality implies that $-v \in \V(P)$. 
For $P_v := \Pi_v(P)$ and $M_v := M / \Z v$ we have
$$\card{\V(P)} \leq \card{\rand P_v \cap M_v} + 4,$$
where equality implies $-v \in \V(P)$. There are now two cases:

\begin{enumerate}
\item Let $w \in \V(P)$ with $w \not= -v$ and $w \not\sim v$. 

Then any lattice point on the boundary of $P$ is in the star set of $v$ 
{\em or} in the star set of $w$ but not away from $w$ {\em or} in $\lin(v,w)$. This implies
$$\card{P_v \cap M_v} + \card{\intr P_w \cap M_w} 
\leq \card{\rand P \cap M} \leq \card{P_v \cap M_v} + \card{\intr P_w \cap M_w} + 2;$$
if the second equality holds, then $-v \in P$. 
\item No such $w$ as in 1. exists. Then:
$$\card{\V(P)} \leq \card{\rand P_v \cap M_v} + 2.$$
\end{enumerate}
\label{casa-bound}
\end{corollary}

Going back to algebraic geometry we derive a generalisation of a theorem 
proven by Casagrande in the nonsingular case \cite[Thm. 2.4]{Cas03a}: 

\begin{corollary}
If $X$ is a $\Q$-factorial Gorenstein toric Fano variety with torus-invariant prime divisor $\Vt$, 
then the Picard numbers satisfy the inequality
\[\rho_X - \rho_{\Vt} \leq 3.\]
\label{pic-bound}
\end{corollary}

Finally using lemmas \ref{fund} and \ref{casa} 
we are now ready to prove the main theorems.

\begin{proof}[Proof of theorems \ref{symmy} and \ref{csymmy}]

Let $P$ be a simplicial reflexive polytope such that there exists a vertex $u \in \V(P^*)$ with $-u \in P^*$. 
Let $F$ be the facet corresponding to $u$ and $F'$ the face defined by $-u$. 
Now define the set $\{v^1, \ldots, v^d\}$ of vertices not in $F$ but in facets 
intersecting $F$ in a codimension two face.  
Lemma \ref{fund} immediately implies that 
$\V(P) \backslash (\V(F) \cup \V(F')) = \{v \in \V(P) \,:\, \pro{u}{v} = 0\} \subseteq \{v^1, \ldots, v^d\}$. 
This yields the bound $\card{\V(P)} \leq 3d$. Now we must consider two special cases:

\begin{enumerate}
\item In order to prove conjecture \ref{mainconj} and thereby finish the proof of theorem \ref{symmy} we may assume 
that $\card{\V(P)} = 3d$ and $d$ is even by remark \ref{enough}. Then \ref{fund} implies that 
$\{x \in \randp \cap M \,:\, \pro{u}{x} = 0\} = \{v^1, \ldots, v^d\}$ is a set of cardinality $d$, and also 
$F'$ must contain $d$ vertices, $F'$ is therefore a $d-1$-dimensional simplex, so $-u \in \V(P^*)$.

Let $e_1, \ldots, e_d$ be the vertices of 
the facet $F$, and $b_1, \ldots, b_d$ the vertices of $F'$ such that $v^1, \ldots, v^d$ are exactly the 
corresponding vertices analogously constructed as in \ref{fund}. Define also $F_i := \conv(v^i, \, 
e_j \,:\, j \not= i)$ for 
$i = 1, \ldots, d$. Then we get the following three facts for $i,j,k \in \{1, \ldots, d\}$:

\smallskip

$\cdot\;${\em Fact 1:} For $i \not=j$: $v^i \sim v^j$ or $v^i + v^j = 0$. 

({\em Proof:} Assume not. Then there exists a $k$ such that $v^i + v^j = v^k \in F_k$. By \ref{prim}{({\rm ii})} 
this implies 
$v^i \in F_k$ or $v^j \in F_k$, giving $v^i = v^k$ or $v^j = v^k$, a contradiction.)

$\cdot\;${\em Fact 2:} For $i$: $e_i + v^i \in \randp$ and $b_i + v^i \in \randp$. 

({\em Proof:} Since $v^i \not\in F_j$ for all $j\not=i$ and $v^i \in F_i$, \ref{fund}(2) yields $e_i \not\sim v^i$. 
By symmetry the same holds for $b_i$.)

$\cdot\;${\em Fact 3:} Let $i,j$ such that $e_i + b_j \in \randp$. 
Then $z(e_i,b_j) = e_i + b_j = v^k$ for some $i \not= k \not=j$. 

({\em Proof:} Since $\pro{u}{e_i + b_j} = 0$, let $e_i + b_j = v^k$ for some $k$. Assume 
$e_i \not\sim v^k$. By \ref{prim} then also $2 e_i + b_j \in \randp$. This implies 
$v^k = 1/2 (b_j + (2 e_i + b_j))$, a contradiction to $v^k \in \V(P)$. 
By symmetry we get $e_i \sim v^k \sim b_j$. By fact 2 necessarily $i \not= k \not=j$.)

\smallskip

Let $i \in \{1, \ldots, d\}$. By fact 2 we can apply lemma \ref{casa} to the vertices $v^i,e_i,b_i$. 
From fact 3 and analysing the possible types in \ref{two} we get that $P \cap \lin(v^i,e_i,b_i)$ must be a 
terminal two-dimensional reflexive polytope, so 
$-e_i = v^i + b_i \in F'$, $-b_i= v^i + e_i \in F$, $v^i = z(-e_i,-b_i) = -e_i + (-b_i)$. 
As this is true for all $i = 1, \ldots, d$, we get $F' = -F$ and $-e_i, -b_i \in \V(P)$. This gives a map
$$s \::\: \{1, \ldots, d\} \to \{1, \ldots, d\}, \;\; i \mapsto s(i), \; \text{ such that } b_{s(i)} :=  -e_i.$$
\begin{enumerate}
\item $s$ is injective, hence a permutation.

\item There are no fixpoints under $s$, i.e., $s(i) \not= i$ for all $\{1, \ldots, d\}$.

\item $-v^i = e_i + b_i \in \rand P$ for all $i \in \{1, \ldots, d\}$.

({\em Proof:} By \ref{casa} it is enough to show that $-v^i \in P$. Assume not.
Fact 1 implies then $v^i \sim v^{s(i)} = z(-e_{s(i)}, -b_{s(i)}) = z(-e_{s(i)}, e_i)$, so by \ref{reid} 
$v^i \sim e_i$, a contradiction to fact 2.)

\item $s \circ s = \id$.

({\em Proof:} Assume there exists an $i \in \{1, \ldots, d\}$ such that for $j := s(i)$ we have $b_{s(j)} \not= b_i$. 
This implies $b_i \sim v^{s(i)} = z(-e_{s(i)}, -b_{s(i)}) = z(b_{s(j)}, e_i)$, so by assumption and 
\ref{reid} $b_i \sim e_i$. This is a contradiction to (c).)
\end{enumerate}

Property (d) implies that $P$ is centrally symmetric. Furthermore 
$s$ is a product of $\frac{d}{2}$ disjoint transpositions in the symmetric group of 
$\{1, \ldots, d\}$. This permutation $s$ 
and the set $\{e_1, \ldots, e_d\}$ of vertices of $F$ uniquely determine $P$, because 
$F' = -F$ and
$v^i = -e_i + e_{s(i)}$ for all $i \in \{1, \ldots, d\}$. 

For any $i \in \{1, \ldots, d\}$ we get $\pro{u}{v^i} = 0$ and 
$\pro{e_i^*}{v^i} = \pro{e_i^*}{-e_i + e_{s(i)}} = -1$. Hence \ref{fund}(3) implies that $e_1, \ldots, e_d$ 
is a $\Z$-basis of $M$. This immediately yields the uniqueness of $P$ up to isomorphism of the lattice.

\item To finish the proof of theorem \ref{csymmy} let $\card{\V(P)} = 3d - 1$, $d$ odd and $P$ be centrally symmetric. 
Then $-u \in \V(P^*)$ and $F' = -F$. By \ref{fund} we can assume that 
$\{x \in \V(P) \cap M \,:\, \pro{u}{x} = 0\} = \{v^1, \ldots, v^{d-1}\}$ is a set of cardinality $d-1$. 
We can again enumerate the vertices of $F$ as $e_1, \ldots, e_d$ such that $v^i$ is a facet of 
$F_i := \conv(v^i,\, e_j \,:\, j \not= i)$ for $i = 1, \ldots, d-1$. There is a map
$$\sigma \::\: \{1, \ldots, d-1\} \to \{1, \ldots, d-1\}, \;\; i \mapsto \sigma(i), \; \text{ such that } 
v^{\sigma(i)} :=  -v^i.$$
$\sigma$ is again just a product of $\frac{d-1}{2}$ disjoint transpositions. Distinguish two cases:

\begin{enumerate}

\item $\pro{u}{v^d} = 0$. 

So there exists exactly one $k \in \{1, \ldots, d-1\}$ such that $v^d = v^k$. We may assume $k = d-1$ and 
$\sigma(d-1) = d-2$. From \ref{fund}(2) and since $v^{d-2} \in \V(P)$ 
we get $y := v^{d-2} + e_{d-2} = z(v^{d-2}, e_{d-2})$ $\in F$, hence $e_{d-2} = y + v^{d-1} = z(y,v^{d-1})$. 
For $j = 1, \ldots, d-3$ we have $e_{d-2} \in F_j$ and $v^{d-1} \not\in F_j$, hence $y \in F_j$ 
by \ref{prim}{({\rm ii})}. 
Since $y$ is away from $e_{d-2}$, this implies $y \in [e_{d-1}, e_d]$. As $e_{d-1} \sim v^d = v^{d-1} 
\sim e_d$ and $v^{d-1} \not\sim y$, it follows that $y \not\in \{e_{d-1}, e_d\}$, hence $y \in ]e_{d-1}, e_d[$. 
Now $e_d \sim v^{d-1} = -y + e_{d-2} = z(-y,e_{d-2})$, so by \ref{reid} also $e_d \sim -y \in ]-e_{d-1}, -e_d[$. 
This implies $e_d \sim (-e_d)$, a contradiction.

\item $\pro{u}{v^d} \not= 0$.

Therefore $v^d \in \V(-F)$. Since $v^d \sim e_i$ for $i = 1, \ldots, d-1$, we have 
$v^d = -e_d$. For $i = 1, \ldots, d-1$ 
lemma \ref{fund}(2) implies again $v^i + e_i \in \randp$, so by central symmetry also 
$v^i - e_{\sigma(i)} \in \randp$. Using 
\ref{casa} we derive $v^i = e_{\sigma(i)} - e_i$. Especially $\pro{u}{v^i} = 0$ and 
$\pro{e_i^*}{v^i} = -1$, hence \ref{fund}(3) implies that $e_1, \ldots, e_d$ is a $\Z$-basis of $M$. This yields the 
uniqueness of $P$, hence $P \cong ([-1,1] \times (S_3)^{\frac{d-1}{2}})^*$.
\end{enumerate}
\end{enumerate}
\end{proof}

\section{Sharp bounds on the number of lattice points}
\label{sec:6}

{\em Throughout the section let $P$ be a $d$-dimensional Fano polytope in $\MR$.}

If $P$ is reflexive, then the number of lattice points in $P$ is 
$h^0(X^*,-K_{X^*})$ for $X^* := X(N, \calN_P)$ 
(see \cite[3.4]{Ful93}). In general there exist very large explicit bounds on $\card{\randp \cap M}$ 
for the class of canonical Fano polytopes (see \cite{ZL91}). However in some cases we can simply count 
integral points modulo $k$, a method that is originally due to Batyrev (see \cite[Lemma 1]{Bat82}).

\begin{definition}{\rm
For $k \in \N$ we have the canonical homomorphism
$$\alpha_k \,:\, M \to M /k M \cong (\Z/ k \Z)^d.$$}
\end{definition}

For a convex set $C \subseteq \MR$ with $C \cap M \not=\emptyset$ 
one easily sees that the minimal $k \in \N_{\geq 1}$ such 
that the restriction of $\alpha_k$ to $C \cap M$ is injective is just the maximal number of lattice points on 
an intersection of $C$ with an affine line. This invariant minus one is called the {\em discrete diameter} of $C$ 
in \cite{Kan98}.

\begin{lemma}
Let $d \geq 2$, $P$ a canonical Fano polytope and 
$B \subseteq \rand P \cap M$ with $\card{[x,y] \cap M} = 2$ for all $x,y \in B$, $x\not=y$, $x \sim y$. 
Let $s$ denote the number of centrally symmetric pairs in $B$. Then
\[\card{B} \leq 2^{d+1} - 2,\;\;\;\;\; s \geq \card{B} + 1 - 2^d.\]
\label{mod}
\end{lemma}

\begin{proof}

We consider the restriction of $\alpha_2$ to $B$. 
As $P$ is canonical, the fibre of $0$ is empty. 
Using the assumption it is also easy to see that the fibre of a non-zero element in $(\Z/2 \Z)^d$ 
has at most two elements, and in the case of equality it consists of one pair of centrally symmetric lattice points in 
$B$. From this the bounds can be derived.

\end{proof}

We immediately get a sharp bound on the number of vertices of a semi-terminal canonical 
Fano polytope (see def. \ref{semi-terminal}). In particular we get a result 
that was proven in the case of a smooth Fano polytope in \cite[Prop. 2.1.11]{Bat99}: 

\begin{corollary}
Let $P$ be a terminal Fano polytope. Then 
$$\card{\randp \cap M} = \card{\VP} \leq 2^{d+1} - 2.$$
If equality holds, then $P$ is centrally symmetric. 
This holds for the terminal reflexive $d$-dimensional standard lattice zonotope $\Zt_d := \conv(\pm [0,1]^d)$, 
see \cite[Proof of thm. 3.21]{DHZ01} for $\Zt^{(d)} = \Zt_{d-1}$.
\label{termpts}
\end{corollary}

The results in \cite{Kas03} show that $\Zt_d$ is even the only terminal Fano polytope with the maximal 
number of vertices for $d \leq 3$. However the computer classification of Kreuzer and Skarke yields 
two non-isomorphic four-dimensional terminal reflexive polytopes with $2^5-2=30$ vertices. 

The second case where counting modulo $k$ works is the class of centrally symmetric reflexive polytopes. 
In the case of a smooth Fano polytope where centrally symmetric pairs of vertices span $\MR$ 
there exists a complete explicit classification that is due to Casagrande (see \cite{Cas03b}). 
However we cannot expect such a result for centrally symmetric reflexive polytopes, 
since by the classification of Kreuzer and Skarke there are $150$ centrally symmetric reflexive polytopes already 
in dimension four. For $d=2$ we have $3$ 
(see \ref{two}) and for $d=3$ there are $13$ (see \cite{Wag95}) 
$d$-dimensional centrally symmetric reflexive polytopes.

There is the conjecture due to Ewald (see \cite{Ewa88}) that {\em any $d$-dimensional smooth Fano polytope 
can be embedded in the unit lattice cube $[-1,1]^d$}. It is proven for $d\leq 4$ by the classification or 
under additional symmetries. It is wrong for simplicial reflexive polytopes, e.g., type $9$ in 
prop. \ref{two} contains $10$ lattice points. 

However {\em for $d \leq 3$ 
we can always embedd a centrally symmetric reflexive polytope $P$ in the unit lattice cube}. 
For this we choose by \ref{klein}(1,2) a $\Z$-basis $b_1, \ldots, b_d$ 
of lattice points in $\rand P^*$, so $P \subseteq \{\sum \lambda_i b^*_i \,:\, \lambda_i \in \{-1,0,1\}\}$ 
for the dual $\Z$-basis $b_1^*, \ldots, b_d^*$.

There is no such result for $d \geq 4$. For instance let $P$ be the $4$-dimensional centrally symmetric reflexive 
polytope $P$ in example \ref{bsp}. Assume $P^*$ could be embedded as a lattice polytope in $[-1,1]^4$. 
Then $([-1,1]^4)^*$ would be a 
lattice subpolytope of $P$ with the same number of vertices. Since $P$ is terminal, this would be an equality, 
a contradiction. This example is taken from \cite{Wir97} where this topic is thoroughly dicussed. In \cite{Wir97} 
there is also a 
characterisation of centrally symmetric reflexive polytopes presented that have the minimal number of $2 d$ vertices, 
and it is shown that these polytopes can always be embedded in the unit lattice cube. 

By such an embedding 
we trivially get that the number of lattice points in the polytope is bounded by $3^d$ with equality only in the case of 
the unit lattice cube. However this is even true in general:

\begin{theorem}
Let $P$ be a centrally symmetric canonical Fano polytope. 
Then
\[\card{P \cap M} \leq 3^d.\]
Any facet of $P$ has at most $3^{d-1}$ lattice points. 
If $P$ is additionally reflexive, then the following statements are equivalent:
\begin{enumerate}
\item $\card{P \cap M} = 3^d$
\item Every facet of $P$ has $3^{d-1}$ lattice points
\item $P \cong [-1,1]^d$ as lattice polytopes
\end{enumerate}
\label{central}
\end{theorem}

\begin{proof}
 
For the bounds we just have to show that $\alpha_3\res{P \cap M}$ is injective. 
So suppose there are $x, y \in P \cap M$ such that $\alpha_3(x) = \alpha_3(y)$. This implies 
$(y-x)/3 \in \intr P \cap M = \{0\}$, hence $x = y$.

The equivalences in the case of a reflexive polytope will be proven in \cite{Nil04}.

\end{proof}


\begin{thebibliography}{XXXXXX}                            

\bibitem[AKMS97]{AKMS97} Avram, A.C.; Kreuzer, M.; Mandelberg, M.; Skarke, H.:  Searching for $K3$ Fibrations. 
Nucl. Phys. B \textbf{494}, 567-589 (1997)

\bibitem[Bat82]{Bat82} Batyrev, V.V.:  Toroidal Fano 3-folds. 
Math. USSR-Izv. \textbf{19}, 13-25 (1982)

\bibitem[Bat91]{Bat91}  Batyrev, V.V.: On the classification of smooth projective toric varieties. 
Tohoku Math. J. \textbf{43}, 569-585 (1991)

\bibitem[Bat94]{Bat94} Batyrev, V.V.: 
Dual polyhedra and mirror symmetry for Calabi-Yau hypersurfaces in 
toric varieties. J. Algebr. Geom. \textbf{3}, 493-535 (1994)

\bibitem[Bat99]{Bat99} Batyrev, V.V.: 
On the classification of toric Fano 4-folds. 
J. Math. Sci., New York \textbf{94}, 1021-1050 (1999)

\bibitem[Bor00]{Bor00} Borisov, A.:  Convex lattice polytopes and cones with few lattice points inside, 
from a birational geometry viewpoint. Preprint, math.AG/0001109 (2000)

\bibitem[Cas03a]{Cas03a} Casagrande, C.:  
Toric Fano varieties and birational morphisms. Int. Math. Res. Not. \textbf{27}, 1473-1505 (2003)

\bibitem[Cas03b]{Cas03b} Casagrande, C.:  Centrally symmetric generators in toric Fano varieties. 
Manuscr. Math. \textbf{111}, No.4, 471-485 (2003)

\bibitem[Con02]{Con02} Conrads, H.:  Weighted projective spaces and reflexive polytopes. 
Manuscripta Math. \textbf{107}, 215-227 (2002)

\bibitem[Deb01]{Deb01} Debarre, O.:  Toric Fano varieties. In: Higher dimensional varieties and rational points, 
lectures of the summer school and conference, Budapest 2001. Bolyai Society Mathematical Studies 12, 
pages 93-132. Berlin: Springer 2001

\bibitem[DHZ01]{DHZ01} Dais, D.I.; Haase, C.; Ziegler, G.M.:  All toric l.c.i.-singularities admit projective 
crepant resolutions. Tohoku Math. J. \textbf{53}, 95-107 (2001)

\bibitem[Ewa88]{Ewa88} Ewald, G.:  On  the  classification  of  toric  
Fano  varieties.  Discrete Comput. Geom. \textbf{3}, 49-54 (1988)

\bibitem[Ewa96]{Ewa96} Ewald, G.:  Combinatorial convexity and algebraic geometry. Graduate texts 
in mathematics 168. New York: Springer 1996

\bibitem[Ful93]{Ful93} Fulton, W.:  
Introduction to toric varieties. Annals of Mathematics Studies 131. Princeton, NJ: 
Princeton University Press 1993

\bibitem[Kan98]{Kan98} Kantor, J.-M.:  Triangulations of integral polytopes and Ehrhart polynomials. 
Contrib. to Alg. and Geom. \textbf{39}, 205-218 (1998)

\bibitem[Kas03]{Kas03} Kasprzyk, A.M.:  Toric Fano 3-folds with terminal singularities. Preprint, math.AG/0311284 (2003)

\bibitem[KS97]{KS97} Kreuzer, M.; Skarke, H.:  On the classification of reflexive polyhedra. 
Commun. Math. Phys. \textbf{185}, 495-508 (1997)

\bibitem[KS98]{KS98} Kreuzer, M.; Skarke, H.:  Classification of reflexive polyhedra in three dimensions. 
Adv. Theor. Math. Phys. \textbf{2}, 853-871 (1998)

\bibitem[KS00a]{KS00a} Kreuzer, M.; Skarke, H.:  Complete classification 
of reflexive polyhedra in four dimensions. 
Adv. Theor. Math. Phys. \textbf{4}, 1209-1230 (2000)

\bibitem[KS00b]{KS00b} Kreuzer, M.; Skarke, H.: Reflexive polyhedra, weights and toric 
Calabi-Yau fibrations. Rev. Math. Phys. \textbf{14}, 343-374 (2002)

\bibitem[KS02]{KS02} Kreuzer, M.; Skarke, H.:  PALP: A package for analyzing lattice polytopes with 
applications to toric geometry. Preprint, math.NA/0204356 (2002)

\bibitem[Nil04]{Nil04} Nill, B.:  Toric Fano varieties with reductive group of automorphisms. 
In preparation (2004)

\bibitem[Oda88]{Oda88} Oda, T.:  Convex bodies and algebraic geometry - An introduction to the theory of toric varieties. 
Ergebnisse der Mathematik und ihrer Grenzgebiete 15. Berlin: Springer 1988

\bibitem[PR00]{PR00} Poonen, B.; Rodriguez-Villegas, F.: Lattice polygons and the number $12$. 
Am. Math. Soc. Monthly \textbf{107}, 238-250 (2000)

\bibitem[Rei83]{Rei83} Reid, M.:  Decomposition of toric morphisms. In: Arithmetic and geometry, Vol. II: Geometry. 
Progress in Mathematics 36, pages 395-418. Boston: Birkh\"auser 1983

\bibitem[Sat00]{Sat00} Sato, H.:  Toward the classification of higher-dimensional Fano varieties. 
Tohoku Math. J., \textbf{52}, 383-413 (2000)

\bibitem[Sti01]{Sti01} Stillwell, J.:  The story of the 120-cell. 
Notices of the Am. Math. Soc. \textbf{48}, 17-24 (2001)

\bibitem[VK85]{VK85} Voskresenskij, V.E.; Klyachko A.A.: 
 Toroidal Fano varieties and root systems. Math. USSR-Izv. 
\textbf{24}, 221-244 (1985) 

\bibitem[Wag95]{Wag95} Wagner, H.:  Gewichtete projektive R\"aume und reflexive Polytope. 
Diplomarbeit (in German). Bochum: Math. Inst. der Ruhr-Universit\"{a}t Bochum 1995

\bibitem[Wir97]{Wir97} Wirth, P.R.:  Zentralsymmetrische reflexive Polytope. Diplomarbeit (in German). 
Bochum: Math. Inst. der 
Ruhr-Universit\"{a}t Bochum 1997

\bibitem[WW82]{WW82} Watanabe, K.; Watanabe, M:  The  classification  
of  Fano  3-folds with torus embeddings. Tokyo J. Math. \textbf{5}, 37-48 (1982)

\bibitem[ZL91]{ZL91} Ziegler, G.M.; Lagarias, J.C.:  Bounds for lattice polytopes containing a fixed 
number of interior points in a sublattice. Can. J. Math. \textbf{43}, 1022-1035 (1991)

\end{thebibliography}
\end{document}